\documentclass[11pt]{amsart}
\usepackage{amscd,amssymb,verbatim}
\theoremstyle{plain}
\newtheorem{Thm}{Theorem}[section]
\newtheorem{Cor}[Thm]{Corollary}
\newtheorem{Lem}[Thm]{Lemma}
\newtheorem{Prop}[Thm]{Proposition}
\newtheorem{Def}[Thm]{Definition}
\newtheorem{remark}[Thm]{Remark}
\newtheorem{notation}[Thm]{Notation}

\numberwithin{equation}{section}

\begin{document}
\title{Weakly null sequences in the Banach space
\(C(K)\)}
\author{I. Gasparis}
\address{Department of Mathematics \\
University of Crete, Knossou Avenue \\
Heracleion 71409, Greece.}
\email{ioagaspa@math.uoc.gr}
\author{E. Odell}
\address{Department of Mathematics \\
The University of Texas at Austin \\
Austin, Texas 78712, U.S.A.}
\email{odell@math.utexas.edu}
\author{B. Wahl}
\address{Department of Mathematics and Computer Science\\
Hanover College \\
P.O. Box 890 \\
Hanover, Indiana 47243, U.S.A.}
\email{wahl@hanover.edu}
\begin{abstract}
The hierarchy of the block bases of transfinite normalized averages of a normalized Schauder
basic sequence is introduced and a criterion
is given for a normalized weakly null sequence in \(C(K)\), the Banach space
of scalar valued functions continuous on the compact metric space \(K\),
to admit a block basis of normalized averages equivalent to the unit vector basis
of \(c_0\), the Banach space of null scalar sequences.
As an application of this criterion, it is shown that
every normalized weakly null sequence in \(C(K)\), for countable \(K\),
admits a block basis of normalized averages equivalent to the unit vector
basis of \(c_0\).
\end{abstract}
\maketitle
\section{Introduction}
We study normalized weakly null sequences
in the spaces \(C(K)\) where \(K\) is a compact metric space.
When \(K\) is uncountable, \(C(K)\) is isomorphic to \(C([0,1])\)
(\cite{mi}, \cite{p}, \cite{be}), while for every countable compact
metric space \(K\) there exist unique countable ordinals \(\alpha\)
and \(\beta\) with \(C(K)\) (linearly) isometric to \(C(\alpha)\) \cite{ms}
and isomorphic (i.e., linearly homeomorphic) to \(C(\omega^{\omega^\beta})\) \cite{bp3}
(in the sequel, for an ordinal \(\alpha\) we let
\(C(\alpha)\) denote \(C([1, \alpha])\), the Banach space of
scalar valued functions, continuous on the ordinal interval
\([1,\alpha]\) endowed with the order topology).

Every normalized weakly null sequence \((f_n)\) in \(C(K)\) for countable \(K\), admits
a basic shrinking subsequence (\cite{bp1}, \cite{dss}) that is, a subsequence \((f_{k_n})\) which is a Schauder basis
for its closed linear span and whose corresponding sequence of biorthogonal
functionals is a Schauder basis for the dual of the closed linear subspace generated by \((f_{k_n})\).

It is shown in \cite{mr} that while \((f_n)\)
must admit an unconditional subsequence in \(C(\omega^\omega)\),
it need not admit an unconditional subsequence in \(C(\omega^{\omega^2})\).

We remark here that if a normalized basic sequence in \(C(K)\) for
countable \(K\) has no weakly null subsequence, then it admits no
unconditional subsequence since such a subsequence would have a
further subsequence {\em equivalent} (this term is explained
below) to the unit vector basis of \(\ell_1\) and \(C(K)\) has
dual isometric to \(\ell_1\) which is separable.

Since \(C(\alpha)\) is \(c_0\)-saturated for all ordinals \(\alpha\) \cite{ps}
(a Banach space is \(c_0\)-saturated provided all of its infinite-dimensional
subspaces contain an isomorph of \(c_0\)), some {\em block basis}
of \((f_n)\) is equivalent to the unit vector basis of \(c_0\).

We recall here that if \((e_n)\) is a Schauder basic sequence in a
Banach space then a non-zero sequence \((u_n)\) is called a block
basis of \((e_n)\), if there exist finite sets \((F_n)\), with
\(\max F_n < \min F_{n+1}\) for all \(n\), and scalars \((a_n)\)
with \(a_i \ne 0 \) for all \(i \in F_n\) and \(n \in \mathbb{N}\)
such that \(u_n = \sum_{i \in F_n } a_i e_i\), for all \(n \in
\mathbb{N}\). We then call \(F_n\) the {\em support} of \(u_n\).
We shall adopt the notation \(u_1 < u_2 < ...\) to indicate that
\((u_n)\) is a block basis of \((e_n)\) such that \(\max \mathrm{
supp } \, u_n < \min \mathrm{ supp } \, u_{n+1}\), for all \(n \in
\mathbb{N}\). We also recall that two basic sequences \((x_n)\),
\((y_n)\) are equivalent provided the map \(T\) sending \(x_n\) to
\(y_n\) for all \(n \in \mathbb{N}\), extends to an isomorphism
between the closed linear spans \(X\) and \(Y\) of \((x_n)\) and
\((y_n)\), respectively. In the case \(T\) only extends to a
bounded linear operator from \(X\) into \(Y\), we say \((x_n)\)
{\em dominates} \((y_n)\).

Our main results are presented mostly in Sections \ref{S1} and \ref{S4}.
We show in Corollary \ref{mC} that if \((f_n)\) is normalized weakly null
in \(C(\omega^{\omega^\xi})\), one can always find \(c_0\)
as a block basis of {\em normalized \(\alpha\)-averages} of
\((f_n)\) for some \(\alpha \leq \xi\), and a quantified
description of \(\alpha\) is given. Note that the proof given in \cite{ps} of
the fact that \(C(\omega^{\omega^\xi})\) is \(c_0\)-saturated
is an existential one that is, it only provides the
existence of a block basis of \((f_n)\) equivalent to the
unit vector basis of \(c_0\) without giving any information
about the support of the blocks or the scalar coefficients
involved. A normalized \(1\)-average of \((f_m)_{m \in M}\)
(where \(M = (m_i)\) is an infinite subsequence of \(\mathbb{N}\))
is a vector \(x = (\sum_{i=1}^{m_1} f_{m_i}) /
\|\sum_{i=1}^{m_1} f_{m_i}\|\). Thus we have that the support
of \(x\) is a maximal \(S_1\)-set in \(M\) where \(S_1\)
is the first Schreier class (see the definition of Schreier classes in the next section).
A \(2\)-average is similarly defined by averaging a block basis
of \(1\)-averages so that the support is a maximal \(S_2\)-set.
This is carried out for all \(\alpha < \omega_1\), as in the construction
of the Schreier classes \(S_\alpha\), yielding the hierarchy of normalized
\(\alpha\)-averages of \((f_n)\). The details are presented in Section \ref{S3}.

Section \ref{S1} includes the following results. We show in
Theorem \ref{C3} and Corollary \ref{T1} that if
a normalized weakly null sequence \((f_n)\)
in \(C(\omega^{\omega^\xi})\) is \(S_\xi\)-unconditional
(see Definition \ref{schr} and the comments after it)
then it admits an unconditional subsequence.
This result, combined with that of \cite{mr} and \cite{o2}
on Schreier unconditional sequences, yields an easier proof
of the aforementioned fact about weakly null sequences
in \(C(\omega^\omega)\) \cite{mr}.
Indeed, as is observed in \cite{mr} (see \cite{o2} for a proof),
every normalized weakly null sequence in a Banach space
admits, for every \(\epsilon > 0\),
a subsequence that is \(S_1\)-unconditional with
constant \(2 + \epsilon\). It follows from this and
Theorem \ref{C3} that every normalized weakly null sequence
in \(C(\omega^\omega)\) admits an unconditional subsequence.
Another consequence of Theorem \ref{C3} is that the example
of a normalized weakly null sequence in \(C(\omega^{\omega^2})\)
without unconditional subsequence \cite{mr}, fails to admit
an \(S_2\)-unconditional subsequence although of course it admits
\(S_1\)-unconditional subsequences. This shows the optimality
of the result in \cite{mr}, \cite{o2} on Schreier unconditional
sequences.

We show in Theorem \ref{C2} that if \((\chi_{G_n})\) is a
weakly null sequence of indicator functions in some space \(C(K)\)
then there exist \(\xi < \omega_1\) and a subsequence of
\((\chi_{G_n})\) which is equivalent to a subsequence of the
unit vector basis of the generalized Schreier space \(X^\xi\)
(\cite{aa}, \cite{ao}) (see Notation \ref{not}).
We thus obtain a quantitative version of Rosenthal's unpublished
result, that a weakly null sequence of indicator
functions in some space \(C(K)\)
admits an unconditional subsequence
(cf. also \cite{amt} and \cite{agr} for another proof of this result).

In Section \ref{S4} we give a sufficient condition for
a normalized weakly null sequence in some \(C(K)\) space
to admit a block basis of normalized averages equivalent
to the unit vector basis of \(c_0\). We show in Theorem \ref{mT}
that if \((f_n)\) is normalized weakly null in \(C(K)\) and
there exist a summable sequence of positive scalars
\((\epsilon_n)\) and a subsequence \((f_{m_n})\)
of \((f_n)\) satisfying
\(\{n \in \mathbb{N} : \, |f_{m_n}(t)| \geq \epsilon_{m_n} \}\)
is finite for all \(t \in K\), then there exist \(\xi < \omega_1\)
and a block basis of normalized \(\xi\)-averages of \((f_n)\)
which is equivalent to the unit vector basis of \(c_0\).
There are two consequences of Theorem \ref{mT}. The first, Corollary
\ref{mC}, has been already discussed. The second one is
Corollary \ref{elton}, which gives
a quantitative
version of a special case of Elton's famous result on extremely
weakly unconditionally convergent sequences \cite{e2}
(cf. also \cite{fo}, \cite{hor}, \cite{an} for related results).
It was shown in \cite{e2} that
if \((x_n)\) is a normalized basic sequence in some
Banach space and the series \(\sum_n |x^*(x_n)|\)
converges for every extreme point \(x^*\) in the
ball of \(X^*\), then some block basis of \((x_n)\)
is equivalent to the unit vector basis of \(c_0\).
We show in Corollary \ref{elton} that if \((f_n)\)
is a normalized basic sequence in some \(C(K)\)
space satisfying \(\sum_n |f_n(t)|\)
converges for all \(t \in K\), then there exist \(\xi < \omega_1\)
and a block basis of normalized \(\xi\)-averages of \((f_n)\)
which is equivalent to the unit vector basis of \(c_0\).

Finally, Sections \ref{S2} and \ref{S3} contain a number
of technical results on \(\alpha\)-averages which are
used in Section \ref{S4}.

Some of the results contained in this paper were obtained
in B. Wahl's thesis \cite{w} written under the supervision of
E. Odell.
\section{Preliminaries}
We shall make use of standard Banach space facts and terminology
as may be found in \cite{lt}. \(c_{00}\) is the vector space
of the ultimately vanishing scalar sequences.
If \(X\) is any set, we let
\([X]^{< \infty}\) denote the set of its finite subsets,
while \([X]\) stands for the set of all infinite subsets
of \(X\). If \(M \in [\mathbb{N}]\), we shall adopt the convenient
notation \(M=(m_i)\) to denote the increasing enumeration of
the elements of \(M\).

A family \(\mathcal{F} \subset [\mathbb{N}]^{< \infty}\) is {\em
hereditary} if \(G \in \mathcal{F}\) whenever \(G \subset F\) and
\(F \in \mathcal{F}\). \(\mathcal{F}\) is {\em spreading} if for
every \(\{m_1 < , \dots , < m_k \} \in \mathcal{F}\) and all
choices \(n_1 < \dots < n_k\) in \(\mathbb{N}\) with \(m_i \leq
n_i\) (\(i \leq k\)), we have that \(\{n_1, \dots , n_k \} \in
\mathcal{F}\). \(\mathcal{F}\) is {\em compact}, if it is compact
with respect to the topology of pointwise convergence in
\([\mathbb{N}]^{< \infty}\). \(\mathcal{F}\) is {\em regular} if
it possesses all three aforementioned properties and contains all
singletons. A regular family \(\mathcal{F}\) is said to be {\em
stable}, provided that \(F \in \mathcal{F}\) is a maximal, under
inclusion, member of \(\mathcal{F}\) if there exists \(n > \max
F\) with \(F \cup \{n\} \notin \mathcal{F}\).

If \(E\) and \(F\) are finite subsets of
\(\mathbb{N}\), we write \(E < F\) when
\(\max E < \min F\).
Given families \(\mathcal{F}_1\) and
\(\mathcal{F}_2\) whose elements are finite subsets
of \(\mathbb{N}\), we define their {\em convolution}
to be the family
\begin{align}
\mathcal{F}_2[\mathcal{F}_1] =
&\{\cup_{i=1}^n G_i: \, n \in \mathbb{N}, \,
G_1 < \dots < G_n, \,
G_i \in \mathcal{F}_1 \, \forall \, i \leq n, \notag \\
&(\min G_i)_{i=1}^n \in \mathcal{F}_2\}
\cup \{\emptyset\}. \notag
\end{align}
It is not hard to see that \(\mathcal{F}_2[\mathcal{F}_1]\)
is regular (resp. stable), whenever each \(\mathcal{F}_i\)
is.

It turns out that for a regular family
\(\mathcal{F}\) there exists a countable
ordinal \(\xi\) such that the \(\xi\)-th Cantor-Bendixson
derivative \(\mathcal{F}^{(\xi)}\) of \(\mathcal{F}\)
is equal to \(\{\emptyset\}\). Hence
\(\mathcal{F}\) is homeomorphic to \([1, \omega^\xi]\), by
the Mazurkiewicz-Sierpinski theorem \cite{ms}.
We then say that \(\mathcal{F}\) is of {\em order}
\(\xi\). If we define
\(\mathcal{F}^{+} =
\{F \in [\mathbb{N}]^{<\infty}: \,
F \setminus \{\min F \} \in \mathcal{F}\}\),
then it is not hard to see, using the
Mazurkiewicz-Sierpinski theorem \cite{ms},
that \( \mathcal{F}^{+}\)
is regular (and stable if \(\mathcal{F}\) is)
of order \(\xi +1\).
It can be shown that if \(\mathcal{F}_i\) is
regular of order \(\xi_i\),
\(i= 1, 2\), then \(\mathcal{F}_2[\mathcal{F}_1]\)
is of order \(\xi_1 \xi_2\).

{\bf Notation}. Given \(\mathcal{F} \subset [\mathbb{N}]^{< \infty}\)
and \(M \in [\mathbb{N}]\), we set
\(\mathcal{F}[M] = \mathcal{F} \cap [M]^{< \infty}\).
Clearly, \(\mathcal{F}[M]\) is hereditary (resp. compact), if
\(\mathcal{F}\) is.

We shall now recall the transfinite definition of
the Schreier families \(S_\xi\), \(\xi < \omega_1\).
First, given a countable ordinal \(\alpha\) we associate
to it a sequence of successor ordinals, \((\alpha_n +1)\),
in the following manner:
If \(\alpha\) is a successor ordinal we let
\(\alpha_n = \alpha - 1\) for all \(n\).
In case \(\alpha\) is a limit ordinal, we choose
\((\alpha_n +1)\) to be a strictly increasing sequence of ordinals
tending to \(\alpha\).

Now set \(S_0 = \{\{n\}: \, n \in \mathbb{N}\} \cup \{\emptyset\}\)
and
\(S_1 = \{F \subset \mathbb{N}: \, |F| \leq \min F\}
\cup \{\emptyset\}\). Note that
\(S_1 = S_1[S_0]\).
Let \(\xi < \omega_1\) and assume \(S_\alpha\) has been defined
for all \(\alpha < \xi\). If \(\xi\) is a successor
ordinal, say
\(\xi = \zeta + 1\), define
\[S_{\xi } = S_1[S_\zeta].\]
In the case \(\xi\) is a limit ordinal,
let \((\xi_n + 1)\) be the
sequence of successor ordinals associated to \(\xi\).
Set
\[S_\xi = \cup_n \{ F \in S_{\xi_n + 1}: \, n \leq \min F\}
\cup \{\emptyset\}.\]
It is shown in \cite{aa} that the Schreier
family \(S_\xi\) is regular of order \(\omega^\xi\)
for all \(\xi < \omega_1\). It is shown in \cite{g}
that the Schreier families are stable.
\begin{Def}[\cite{mr}, \cite{o2}] \label{schr}
A normalized
basic sequence \((x_n)\) in a Banach space is said to be
Schreier unconditional,
if there exists a constant \(C>0\)
such that \(\|\sum_{n \in F} a_n x_n\| \leq C \|\sum_n a_n x_n\|\),
for every \(F \subset \mathbb{N}\) with \(|F| \leq \min F\),
and all choices of finitely supported scalar sequences \((a_n)\).
\end{Def}
It has been already mentioned in the introductory section that
every normalized weakly null sequence admits, for every \(\epsilon
> 0\), a subsequence that is Schreier unconditional with constant
\(2 + \epsilon\).

The concept of Schreier unconditionality
can be generalized in the following manner: Consider
a hereditary family \(\mathcal{F}\) of finite subsets of \(\mathbb{N}\)
containing the singletons.
A normalized basic sequence \((x_n)\) is now called
\(\mathcal{F}\)-{\em unconditional}, if
there exists a constant \(C>0\)
such that \(\|\sum_{n \in F} a_n x_n\| \leq C \|\sum_n a_n x_n\|\),
for every \(F \in \mathcal{F} \)
and all choices of finitely supported scalar sequences \((a_n)\).
\section{Upper Schreier estimates} \label{S1}
In this section we show that every normalized weakly null sequence
in \(C(K)\), \(K\) a countable compact metric space, admits a
subsequence dominated by a subsequence of the unit vector basis
of a certain Schreier space (see the relevant definition after
the statement of Theorem \ref{C1}).

Recall, \cite{bp3}, that for every
countable compact metric space \(K\), there exists a unique
countable ordinal \(\alpha\) with \(C(K)\) isomorphic to
\(C(\omega^{\omega^\alpha})\).
Since most of the properties
of weakly null sequences in \(C(K)\) that we shall be interested in,
are isomorphic invariants, there will be no loss of generality
in assuming that \(K = [1, \omega^\xi]\), for some \(\xi < \omega_1\).
As is has been already mentioned in the previous
section, every regular family \(\mathcal{F}\) of order \(\xi \)
(this means \(\mathcal{F}^{(\xi)} = \{\emptyset\}\)) is homeomorphic
to the ordinal interval \([1, \omega^\xi]\). Moreover, it is easy
to construct by transfinite induction, a regular family of order \(\xi\),
for all \(\xi < \omega_1\). We can thus identify \(C(\omega^\xi)\)
with \(C(\mathcal{F})\), for every regular family of order \(\xi\).

The advantage of such a representation is that one can easily construct
a monotone, shrinking Schauder basis of \(C(\mathcal{F})\), the
so-called {\em node basis} \cite{ajo}. Indeed, let \((\alpha_n)_{n=1}^\infty\)
be an enumeration of the elements of \(\mathcal{F}\), compatible
with the natural partial ordering of \(\mathcal{F}\) given by initial segment
inclusion. This means that whenever \(\alpha_m\) is a proper
initial segment of \(\alpha_n\), then \(m < n\). In particular,
\(\alpha_1 = \emptyset\). Such an enumeration is for instance,
the anti-lexicographic enumeration of the elements of \(\mathcal{F}\),
i.e., \(F \prec G\) if and only if either \(\max F < \max G\), or
\(F \setminus \{\max F\} \prec G \setminus \{\max G\}\),
for all \(F\), \(G\) in \(\mathcal{F}\).

Given \(\alpha \in \mathcal{F}\), set
\(G_\alpha = \{ \beta \in \mathcal{F} : \, \alpha \leqslant \beta\}\),
where \(\alpha \leqslant \beta\) means that \(\alpha\) is an
initial segment of \(\beta\). Clearly, \(G_\alpha\) is a
clopen subset of \(\mathcal{F}\) for every
\(\alpha \in \mathcal{F}\).
The sequence \((\chi_{G_{\alpha_n}})_{n=1}^\infty\)
is called the node basis of \(C(\mathcal{F})\).
It is not hard to check that
\((\chi_{G_{\alpha_n}})_{n=1}^\infty\) is
a normalized, monotone, shrinking Schauder basis for
\(C(\mathcal{F})\) \cite{ajo}.
\begin{Prop} \label{node}
Let \(\mathcal{F}\) be a regular family and
\(u_1 < u_2 < \dots\) be a block basis
of the node basis \((\chi_{G_{\alpha_n}})_{n=1}^\infty\)
of \(C(\mathcal{F})\). Then there exist positive integers
\(n_1 < n_2 < \dots\) with the following property:
For every \(\gamma \in \mathcal{F}\),
\(\{n_i : \, i \in \mathbb{N}, \, u_{n_i}(\gamma) \ne 0\} \in \mathcal{F}^{+}\).
\end{Prop}
\begin{proof}
Define \(F_n = \{ \alpha_i : \, i \in \mathrm{ supp } \, u_n\}\),
for all \(n \in \mathbb{N}\). Clearly, the \(F_n\)'s are pairwise
disjoint, finite subsets of \(\mathcal{F}\). We observe that
whenever \(\alpha_i \in F_n\) and \(\alpha_j \in F_m\) satisfy
\(\alpha_i \leqslant \alpha_j\), then \(n \leq m\). This is so
since \(\alpha_i \leqslant \alpha_j\) implies that \(i \leq j\)
and, subsequently, that \(u_n \leq u_m\). Hence, \(n \leq m\).

We next choose inductively, integers \(2 = n_1 < n_2 < \dots\)
such that \(\max \beta < n_{i+1} \) for every \( \beta \in
F_{n_i}\) and all \(i \in \mathbb{N}\) (where, \(\max \beta\)
denotes the largest element of the finite subset \(\beta\) of
\(\mathbb{N}\)). We claim \((n_i)\) is as desired. Indeed, let
\(\gamma \in \mathcal{F}\). Then
\[\{n_i : \, i \in \mathbb{N}, \, u_{n_i}(\gamma) \ne 0\} \subset
\{n_i : \, i \in \mathbb{N}, \, \exists \, \beta \in F_{n_i}, \,
\beta \leqslant \gamma\},\] for writing \(u_{n_i} = \sum_{\beta
\in F_{n_i}} \lambda_\beta \chi_{G_\beta}\) for some suitable
choice of scalars \((\lambda_\beta)_{ \beta \in F_{n_i}}\), we see
that \(u_{n_i}(\gamma) \ne 0\) implies \(\chi_{G_\beta}(\gamma) =
1\), for some \(\beta \in F_{n_i}\) with \(\beta \leqslant
\gamma\). In particular, \(\{n_i : \, i \in \mathbb{N}, \,
u_{n_i}(\gamma) \ne 0\}\) is finite. Let now \(\{n_{i_1} < ,
\dots, < n_{i_k}\}\) be an enumeration of \(\{n_i : \, i \in
\mathbb{N}, \, u_{n_i}(\gamma) \ne 0\}\), and choose \(\beta_j \in
F_{n_{i_j}}\) with \(\beta_j \leqslant \gamma\), for all \(j \leq
k\). Since \(\{\beta_1, \dots , \beta_k\}\) is well-ordered with
respect to the partial ordering \(\leqslant\) of \(\mathcal{F}\)
(all the \(\beta_j\)'s are initial segments of \(\gamma\)), our
preliminary observation yields \(\beta_1 <  \dots < \beta_k\).
Note that \(\beta_1 \ne \emptyset\). By the choices made, \(\max
\cup \, \beta_j < n_{i_j + 1} \leq n_{i_{j+1}}\) for all \(j \leq
k\). Because \(\mathcal{F}\) is hereditary and spreading, we infer
that \(\{n_{i_2}, \dots , n_{i_k} \} \in \mathcal{F}\) whence
\(\{n_i : \, i \in \mathbb{N}, \, u_{n_i}(\gamma) \ne 0\} \in
\mathcal{F}^{+}\), as required.
\end{proof}
\begin{Cor} \label{C0}
Suppose \(K\) is homeomorphic to
\([1, \omega^\xi]\), \(\xi < \omega_1\), and that
\((f_i)\) is a normalized weakly null sequence in \(C(K)\).
Let \(\mathcal{F}\) be a regular family of order \(\xi\).
Then for every \(N \in [\mathbb{N}]\) and every non-increasing
sequence of positive scalars \((\epsilon_i)\),
there exists \(M \in [N]\), \(M= (m_i)\),
such that
for every \(t \in K\) the set
\(\{m_i : \, i \in \mathbb{N}, \, |f_{m_i}(t)| \geq \epsilon_i \}\)
belongs to \(\mathcal{F}^{+}\).
\end{Cor}
\begin{proof}
We identify \(C(\mathcal{F})\) with \(C(K)\)
and apply Proposition \ref{node}
to find a normalized, shrinking, monotone
Schauder basis \((e_i)\) for \(C(K)\)
with the following property: For every block basis
\(u_1 < u_2 < \dots \) of \((e_i)\) there exist
positive integers \(n_1 < n_2 < \dots\)
such that for all \(t \in K\),
\(\{n_i : \, i \in \mathbb{N}, \, u_{n_i}(t) \ne 0 \}
\in \mathcal{F}^{+}\).

Now let \((f_i)\) be normalized weakly null in \(C(K)\).
A classical perturbation result \cite{bp1} yields a subsequence
\((f_{l_i})\) of \((f_i)\) and a block basis \((u_i)\) of
\((e_i)\), \(u_1 < u_2 < \dots\), such that \(l_i \in N\) and
\(\| f_{l_i} - u_i \| < \epsilon_i /2\), for all \(i \in \mathbb{N}\).
We next choose positive integers \(n_1 < n_2 < \dots\)
such that \(\{n_i : \, i \in \mathbb{N}, \, u_{n_i}(t) \ne 0 \}
\in \mathcal{F}^{+}\), for all \(t \in K\). Set
\(m_i = l_{n_i}\), for all \(i \in \mathbb{N}\).
It is not hard to check using the spreading property of \(\mathcal{F}\),
that \(M = (m_i)\) satisfies the desired conclusion.
\end{proof}
\begin{notation} \label{not}
Let \(\mathcal{F}\) be
a regular family and let
\((e_i)\) denote the unit vector basis of \(c_{00}\).
We define a norm \(\| \cdot \|_{\mathcal{F}}\) on
\(c_{00}\) by the rule
\[\bigl \|\sum_i a_i e_i \bigr \|_{\mathcal{F}} = \sup
\bigl \{ \sum_{i \in F} |a_i| : \, F \in \mathcal{F}\bigr \},
\text{ for all } (a_i) \in c_{00}.\]
The completion of \((c_{00}, \| \cdot \|_{\mathcal{F}})\)
is a Banach space having \((e_i)\) as a normalized, unconditional,
shrinking, monotone Schauder basis (see \cite{aa}, \cite{ao}).
When \(\mathcal{F} = S_\xi\), the \(\xi\)-th Schreier class,
we obtain the generalized Schreier space \(X^\xi\) introduced
in \cite{aa}, \cite{ao}.
\end{notation}
Our next result yields that every
normalized weakly null sequence in
\(C(\omega^{\omega^\xi})\) admits a subsequence
dominated by a subsequence of the unit vector basis
of the generalized Schreier space \(X^\xi\).
\begin{Prop} \label{P1}
Suppose \(K\) is homeomorphic to
\([1, \omega^\xi]\), \(\xi < \omega_1\), and that
\((f_i)\) is a normalized weakly null sequence in \(C(K)\).
Let \(\mathcal{F}\) be a regular family of order \(\xi\).
Given \(0 < \epsilon < 1\), there exists \(M \in [\mathbb{N}]\),
\(M = (m_i)\), such that
\begin{align}
\bigl \|\sum_i a_i f_{m_i} \bigr \| &\leq \frac{2}{1-\epsilon}
\sup \bigl \{ \bigl \| \sum_{i \in F} a_i f_{m_i} \bigr \| : \, F \subset \mathbb{N}, \,
(m_i)_{i \in F} \in \mathcal{F} \bigr \} \notag \\
&\leq \frac{2}{1-\epsilon} \bigl \|\sum_i a_i e_{m_i} \bigr \|_{\mathcal{F}},
\text{ for all } (a_i) \in c_{00}. \notag
\end{align}
\end{Prop}
\begin{proof}
We may assume that \((f_i)\) is \(2\)-basic.
Choose a decreasing sequence of positive scalars \((\epsilon_i)\)
such that \(\sum_i \epsilon_i < \epsilon / 3 \).
We next choose \(M \in [\mathbb{N}]\), \(M= (m_i)\),
satisfying the conclusion of Corollary \ref{C0}
applied to \((f_i)\) and the scalar sequence \((\epsilon_i)\).

Let \((a_i) \in c_{00}\) be such that
\(\|\sum_i a_i f_{m_i} \| =1\), and let \(t \in K\)
satisfy \(|\sum_i a_i f_{m_i} (t) | =1\).
Since \(\{m_i : \, i \in \mathbb{N}, \, |f_{m_i}(t) | \geq \epsilon_i \}\)
belongs to \(\mathcal{F}^{+}\), we obtain
\[1 \leq 2 \sup \bigl \{ \bigl \| \sum_{i \in F} a_i f_{m_i} \bigr \| : \, F \subset \mathbb{N}, \,
(m_i)_{i \in F} \in \mathcal{F} \bigr \} + \epsilon\]
from which the assertion of the proposition follows.
\end{proof}
\begin{remark}
S. Argyros has discovered an alternate proof
of Corollary \ref{C0}. He shows that given
a weakly null sequence \((f_i)\) in \(C(\omega^\xi)\)
and a summable sequence
of positive scalars \((\epsilon_i)\) then, by identifying
\(C(\omega^\xi)\) with \(C(\mathcal{F})\), one can
select positive integers \(1 = m_1 < m_2 < \dots\)
such that if \(|f_{m_i}(F)| \geq \epsilon_i\)
for some \(i \in \mathbb{N}\) and \(F \in \mathcal{F}\),
then \(F \cap (m_{i-1}, m_{i+1}) \ne \emptyset\)
(\(m_0 = 0\)). Therefore, \(\{m_{2i} : \, i \in \mathbb{N}, \, |f_{m_{2i}} (F)| \geq \epsilon_{2i} \}
\in \mathcal{F}^{+}\), for every \(F \in \mathcal{F}\)
which clearly implies Corollary \ref{C0}.
\end{remark}
\begin{remark}
Proposition 9 and Lemma 13 in \cite{l} yield
that for a normalized weakly null sequence \((f_i)\)
in \(C(\omega^\xi)\) there exist a subsequence
\((f_{m_i})\), a compact hereditary family
\(\mathcal{D}\) with \(\mathcal{D}^{(\xi +1)} = \emptyset\)
and a constant \(d > 0\) such that
\(\|\sum_i a_i f_{m_i}\| \leq d
\sup \{\|\sum_{i \in A} a_i f_{m_i}\|: \, A \in \mathcal{D}\}\)
for every \((a_i) \in c_{00}\).
\end{remark}
\begin{Thm} \label{C3}
Suppose \(K\) is homeomorphic to
\([1, \omega^\xi]\), \(\xi < \omega_1\), and that
\((f_i)\) is a normalized weakly null sequence in \(C(K)\).
Let \(\mathcal{F}\) be a regular family of order \(\xi\).
Assume \((f_i)\) is \(\mathcal{F}\)-unconditional.
Then \((f_i)\) has an unconditional subsequence.
\end{Thm}
\begin{proof}
Suppose \((f_i)\) is \(\mathcal{F}\)-unconditional with constant
\(C > 0\). This means that \(\|\sum_{i \in F} a_i f_i \| \leq C
\|\sum_i a_i f_i \|\), for all \(F \in \mathcal{F}\) and every
\((a_i) \in c_{00}\). Let \(M= (m_i)\) satisfy the conclusion of
Proposition \ref{P1}, for \((f_i)\) and \(\mathcal{F}\) with
\(\epsilon = 1/2\). We claim that \((f_{m_i})\) is unconditional.
Indeed, let \((a_i) \in c_{00}\) and \(I \in [\mathbb{N}]\).
Proposition \ref{P1} yields
\[\bigl \|\sum_{i \in I} a_i f_{m_i} \bigr \| \leq
4 \sup \bigl \{ \bigl \|\sum_{i \in F \cap I} a_i f_{m_i} \bigr \| : \, F \subset \mathbb{N}, \,
(m_i)_{i \in F} \in \mathcal{F} \bigr \}.\]
Since \(\mathcal{F}\) is hereditary and \((f_i)\) is \(\mathcal{F}\)-unconditional,
we have that
\[ \bigl \|\sum_{i \in F \cap I} a_i f_{m_i} \bigr \| \leq C
\bigl \|\sum_i a_i f_{m_i} \bigr \|,
\text{ whenever } (m_i)_{i \in F} \in \mathcal{F}.\]
Therefore, \(\|\sum_{i \in I} a_i f_{m_i} \| \leq 4C \|\sum_i a_i f_{m_i}\|\)
which proves the claim.
This completes the proof.
\end{proof}
From Theorem \ref{C3} we easily obtain the next
\begin{Cor} \label{T1}
A normalized weakly null sequence in
\(C(\omega^{\omega^\xi})\), \(\xi < \omega_1\), admits
an unconditional subsequence if, and only if, it admits
a subsequence which is \(S_\xi\)-unconditional.
\end{Cor}
\begin{Thm} \label{C1}
Let \((f_i)\) be a normalized weakly null sequence in
\(C(\omega^{\omega^\xi})\), \(\xi < \omega_1\).
Assume that \((f_i)\) is an \(\ell_1^\xi\)-spreading model.
Then \((f_i)\) admits a subsequence equivalent to a subsequence
of the unit vector basis of \(X^\xi\), the generalized Schreier
space of order \(\xi\) (see Notation \ref{not}).
\end{Thm}
We recall that
\(X^0 = c_0\) while \(X^1\) was implicitly considered by
Schreier \cite{sc}. The generalized Schreier spaces
\(X^\xi\), \(\xi < \omega_1\), were introduced in
\cite{aa}, \cite{ao}.
They can be thought as the
the higher ordinal unconditional analogs of \(c_0\).

We also recall (\cite{amt}), that a normalized basic sequence
\((x_n)\) is said to be an \(\ell_1^\xi\)-{\em spreading model},
\(\xi < \omega_1\), if there is a constant \(\delta > 0\) such
that \(\|\sum_{n \in F} a_n x_n\| \geq \delta \sum_{n \in F}
|a_n|\), for every \(F \in S_\xi\) and all choices of scalars
\((a_n)_{n \in F}\). Saying \((x_n)\) is an \(\ell_1^1\)-spreading
model means that \(\ell_1\) is a spreading model for the space
generated by some subsequence of \((x_n)\), in the sense of
\cite{bs}, \cite{bl}, \cite{o1}. \(\ell_1^\xi\)-spreading models
are instrumental in the study of {\em asymptotic}
\(\ell_1\)-spaces \cite{otw}. It is shown in \cite{ag} that a
weakly null sequence which is an \(\ell_1^\xi\)-spreading model,
admits a subsequence which is \(S_\xi\)-unconditional. The unit
vector basis of \(X^\xi\) is an \(\ell_1^\xi\)-spreading model
with constant \(\delta = 1\).
\begin{proof}[Proof of Theorem \ref{C1}]
We first apply Proposition \ref{P1} with \(\epsilon = 1/2\),
to obtain an infinite subset \(M=(m_i)\)
of \(\mathbb{N}\)
with
\(\|\sum_i a_i f_{m_i} \| \leq 4
\|\sum_i a_i e_{m_i} \|_{S_\xi}\)
for all \((a_i) \in c_{00}\),
where \((e_i)\) denotes the unit vector basis of
\(X^\xi\).
On the other hand, as \((f_i)\) is an
\(\ell_1^\xi\)-spreading model, there exists
a constant \(\delta > 0\) such that
\[ \|\sum_i a_i f_{m_i} \| \geq \delta
\|\sum_i a_i e_{m_i} \|_\xi,
\text{ for all }
(a_i) \in c_{00}.\]
We infer from the preceding inequalities
that \((f_{m_i})\) and \((e_{m_i})\)
are equivalent.
\end{proof}
Our final result in this section yields a quantitative version of
Rosenthal's result, that a weakly null (in \(C(K)\)) sequence of
indicator functions of clopen subsets of a compact Hausdorff space
\(K\), admits an unconditional subsequence (cf. also \cite{amt}
and \cite{agr} for another proof of this result).
\begin{Thm} \label{C2}
Let \(K\) be a compact Hausdorff space.
Suppose that \((f_n)\) is a normalized weakly null sequence
in \(C(K)\) such that there exists \(\epsilon > 0\)
with the property \(f_n (t) = 0\) or \(|f_n(t)| \geq \epsilon\)
for all \(t \in K\) and \(n \in \mathbb{N}\).
Then there exist \(\xi < \omega_1\) and
a subsequence of \((f_n)\) equivalent to a
subsequence of the natural Schauder basis
of \(X^\xi\).
\end{Thm}
\begin{proof}
We first employ the results of \cite{aa}
in order to find the smallest countable ordinal
\(\eta\) for which there is a subsequence \((f_{m_n})\)
of \((f_n)\),
such that no subsequence of \((f_{m_n})\)
is an \(\ell_1^\eta\)-spreading
model. Such an ordinal exists because
\((f_n)\) is weakly null.
We claim that \(\eta\) is a successor ordinal.
To see this we shall need a result from \cite{ag}
(Corollary 3.6)
which states that a weakly null sequence \((f_n)\)
in a \(C(K)\) space admits a subsequence which is
an \(\ell_1^\alpha\)-spreading model, for some
\(\alpha < \omega_1\) if, and only if, there exist
a constant \(\delta > 0\) and \(L \in [\mathbb{N}]\),
\(L=(l_n)\), so that for every \(F \in S_\alpha\)
there exists \(t \in K\) satisfying
\(|f_{l_n}(t)| \geq \delta\), for all \(n \in F\).

Define \(G_n = \{t \in K : \, f_n(t) \ne 0\}\).
Our assumptions yield
\(G_n = \{t \in K : \, |f_n(t)| \geq \epsilon \}\),
for all \(n \in \mathbb{N}\).
Observe that for every
\(\alpha < \eta\) and \(P \in [\mathbb{N}]\),
there exists \(Q \in [P]\),
\(Q=(q_n)\), so that \((f_{q_n})\) is
an \( \ell_1^\alpha\)-spreading model.
It follows now, from the previously cited
result of \cite{ag}, that
for every
\(\alpha < \eta\) and \(P \in [\mathbb{N}]\),
there exists \(Q \in [P]\),
\(Q=(q_n)\), so that
for every \(F \in S_\alpha\),
\(\cap_{n \in F} G_{q_n} \ne \emptyset\).
This in turn
yields that every subsequence of \((f_{m_n})\)
admits, for every \(\alpha < \eta\), a further subsequence
which is an \( \ell_1^\alpha\)-spreading model with constant
independent of \(\alpha\) and the particular subsequence.
Were \(\eta\) a limit ordinal, we would have that some
subsequence of \((f_{m_n})\)
is an \(\ell_1^\eta\)-spreading model, contrary
to our assumption.

Hence, \(\eta = \xi +1\), for some \(\xi < \omega_1\).
Let \((e_n)\) be the natural basis of \(X^\xi\). We
show that some subsequence of \((f_{m_n})\)
is equivalent to a subsequence of \((e_n)\).
Because \(\xi < \eta\),
we can assume without loss of generality, after
passing to a subsequence if necessary, that
\((f_{m_n})\) is an \(\ell_1^\xi\)-spreading
model and thus there exists a constant \(\rho > 0\)
such that
\(\|\sum_n a_n f_{m_n}\| \geq \rho
\|\sum_n a_n e_n \|_{S_\xi}\)
for all \((a_n) \in c_{00}\).
Define
\[ \mathcal{F} = \{F \in [\mathbb{N}]^{< \infty}: \,
\cap_{i \in F} G_{m_i} \ne \emptyset\}.\]
Clearly, \(\mathcal{F}\) is hereditary. It is shown in
\cite{ag},
based on the fact that no subsequence of \((f_{m_n})\)
is an \(\ell_1^{\xi + 1}\)-spreading model,
that there exist \(L \in [\mathbb{N}]\), \(L=(l_n)\),
and \(d \in \mathbb{N}\) so that every member of
\(\mathcal{F}[L]\) is contained in the union of
\(d\) members of \(S_\xi[L]\).
Let \(k_n = m_{l_n}\), for all \(n \in \mathbb{N}\).
We deduce from our preceding work that
\(\|\sum_n a_n f_{k_n} \| \leq d
\|\sum_n a_n e_{l_n}\|_{S_\xi}\), for every
\((a_n) \in c_{00}\).
Therefore, \((f_{k_n})\) and \((e_{l_n})\)
are equivalent.
\end{proof}
\section{Normalized averages of a basic sequence} \label{S2}
Let \(\vec{s} = (e_n)\) be a normalized basic sequence
in a Banach space, and let \(\mathcal{F}\)
be a regular and stable family.
We shall introduce an hierarchy
\(\{(\alpha_n^{\mathcal{F}, \vec{s}, M})_{n=1}^\infty, \,
M \in [\mathbb{N}], \, \alpha < \omega_1 \}\)
of normalized
block bases of \(\vec{s}\), similar to that
of the repeated averages introduced in \cite{amt}.
The latter however consists of convex block bases of
\(\vec{s}\), not necessarily normalized.

We fix a normalized basic sequence \(\vec{s}= (e_n)\) and
a regular and stable family \(\mathcal{F}\).
To simplify our notation, we shall write
\(\alpha_n^M\) instead of \(\alpha_n^{\mathcal{F}, \vec{s}, M}\).
We shall next define, by transfinite induction on \(\alpha < \omega_1\),
a family of normalized block bases
\((\alpha_n^M)_{n=1}^\infty\) of \(\vec{s}\),
where \(M \in [\mathbb{N}]\), so that the following properties
are fulfilled for every
\(\alpha < \omega_1\) and \(M \in [\mathbb{N}]\):
\begin{enumerate}
\item \(\alpha_n^M < \alpha_{n+1}^M\), for all \(n \in \mathbb{N}\).
\item \(M = \cup_n \mathrm{ supp } \, \alpha_n^M\),
for all \(M \in [\mathbb{N}]\).
\end{enumerate}
If \(\alpha = 0\) and \(M= (m_n)\) set
\(\alpha_n^M = e_{m_n}\), for all \(n \in \mathbb{N}\).

Suppose \((\beta_n^N)_{n=1}^\infty\) has been defined
so that \((1)\) and \((2)\), above, are satisfied for
all \(\beta < \alpha\) and \(N \in [\mathbb{N}]\).
Let \(M \in [\mathbb{N}]\). In order to define
\((\alpha_n^M)_{n=1}^\infty\), assume first that
\(\alpha\) is successor, say \(\alpha = \beta +1\).
Let \(k_1\) be the unique integer such that the set
\(\{\min \mathrm{ supp } \, \beta_i^M : \, i \leq k_1\}\)
is a maximal member of \(\mathcal{F}\). We define
\[ \alpha_1^M = \bigl (\sum_{i=1}^{k_1} \beta_i^M \bigr )
\, \bigl / \,
\bigl \|\sum_{i=1}^{k_1} \beta_i^M \bigr \|.\]
Suppose that \(\alpha_1^M < \dots < \alpha_n^M\)
have been defined and that the union of their supports
forms an initial segment of \(M\). Set
\[M_{n+1} = \{ m \in M : \,
\max \mathrm{ supp } \, \alpha_n^M < m\}.\]
Let \(k_{n+1}\) be the unique integer such that the set
\(\{\min \mathrm{ supp } \, \beta_i^{M_{n+1}} : \, i \leq k_{n+1}\}\)
is a maximal member of \(\mathcal{F}\). We define
\[ \alpha_{n+1}^M = \bigl (\sum_{i=1}^{k_{n+1}} \beta_i^{M_{n+1}}
\bigr ) \, \bigl / \, \bigl \|\sum_{i=1}^{k_{n+1}}
\beta_i^{M_{n+1}} \bigr \|.\] This completes the definition of
\((\alpha_n^M)_{n=1}^\infty\) when \(\alpha\) is a successor
ordinal. Note that the construction described above can be carried
out because \(\mathcal{F}\) is stable. \((1)\) and \((2)\) are now
satisfied by \((\alpha_n^M)_{n=1}^\infty\).

Now suppose \(\alpha\) is a limit ordinal.
Let \((\alpha_n + 1)\) be the sequence of successor ordinals
associated to \(\alpha\). Let \(M \in [\mathbb{N}]\)
and set \(m_1 = \min M\).
In case \(m_1 = 1\), set \(\alpha_1^M = e_1\).
If \(m_1 > 1\),
define
\[\alpha_1^M = u^M \, / \, \|u^M\|, \text{ where }
u^M = (1/m_1) e_{m_1} + [\alpha_{m_1}]_1^{M \setminus \{m_1\}} .\]
Suppose that \(\alpha_1^M < \dots < \alpha_n^M\)
have been defined and that the union of their supports
forms an initial segment of \(M\). Set
\[M_{n+1} = \{ m \in M : \,
\max \mathrm{ supp } \, \alpha_n^M < m\}\]
and \(m_{n+1} = \min M_{n+1}\). Define
\begin{align}
&\alpha_{n+1}^M = u^{M_{n+1}} \, / \, \|u^{M_{n+1}}\|,
\text{ where } \notag \\
&u^{M_{n+1}} = (1/m_{n+1}) e_{m_{n+1}} +
[\alpha_{m_{n+1}}]_1^{M_{n+1} \setminus \{m_{n+1} \}}. \notag
\end{align}
Note that \(\alpha_{n+1}^M = \alpha_1^{M_{n+1}}\).
This completes the definition \((\alpha_n^M)_{n=1}^\infty\)
when \(\alpha\) is a limit ordinal. It is clear that
\((1)\) and \((2)\) are satisfied.
\begin{remark}
In case \(\mathcal{F} = S_1\), the first Schreier family,
it is not hard to see that \(\mathrm{ supp } \, \alpha_n^M \in S_\alpha\),
for all \(\alpha < \omega_1\), all \(M \in [\mathbb{N}]\) and all
\(n \in \mathbb{N}\).
\end{remark}
The next lemma is an immediate consequence of the preceding
definition.
\begin{Lem} \label{L11}
Let \(\alpha < \omega_1\), \(M \in [\mathbb{N}]\) and
\(n \in \mathbb{N}\). Then there exists \(N \in [\mathbb{N}]\)
such that \(\alpha_n^M = \alpha_1^N\).
\end{Lem}
Our next result will be applied later,
in conjunction with the infinite Ramsey theorem,
in order to determine if there exists a block basis
of the form \((\alpha_n^M)\), equivalent to the \(c_0\)-basis.
\begin{Lem} \label{L12}
Let \(\alpha < \omega_1\), \(M \in [\mathbb{N}]\)
and \(n \in \mathbb{N}\).
Let \(L_i \in [\mathbb{N}]\) and \(k_i \in \mathbb{N}\),
for \(i \leq n\), be so that
\(\alpha_{k_1}^{L_1} < \dots < \alpha_{k_n}^{L_n}\)
and
\(\cup_{i=1}^n \mathrm{ supp } \, \alpha_{k_i}^{L_i}\)
is an initial segment of \(M\).
Then \(\alpha_i^M = \alpha_{k_i}^{L_i}\),
for all \(i \leq n\).
\end{Lem}
\begin{proof}
By Lemma \ref{L11} we may assume that \(k_i = 1\)
for all \(i \leq n\). We prove the assertion of the lemma
by transfinite induction on \(\alpha\).
The case \(\alpha = 0\) is trivial.
Suppose the assertion holds for all ordinals smaller
than \(\alpha\), and all \(M \in [\mathbb{N}]\)
and \(n \in \mathbb{N}\). Let \(M \in [\mathbb{N}]\).
We prove the assertion for \(\alpha\) by induction on
\(n\). If \(n=1\), we first consider the case of
\(\alpha\) being a successor ordinal, say
\(\alpha = \beta +1\).
We know from the definitions that
\[ \mathrm{ supp } \, \alpha_1^{L_1} = \cup_{j=1}^{p_1}
\mathrm{ supp } \, \beta_j^{L_1},\]
where
\(\{\min \mathrm{ supp } \, \beta_j^{L_1} : \, j \leq p_1\}\)
is a maximal member of \(\mathcal{F}\).
In particular, the set \(\cup_{j=1}^{p_1}
\mathrm{ supp } \, \beta_j^{L_1}\)
is an initial segment of \(M\).
The induction hypothesis on \(\beta\) now implies that
\(\beta_j^M = \beta_j^{L_1}\), for all \(j \leq p_1\).
It follows now that \(\alpha_1^M = \alpha_1^{L_1}\).

To complete the case \(n=1\), we
consider the possibility that \(\alpha\) is a limit ordinal.
Let \((\alpha_n +1)\) be the sequence of ordinals associated
to \(\alpha\) and suppose that \(m = \min M\). Then
\(m = \min \mathrm{ supp } \, \alpha_1^{L_1}\)
and so \(m = \min L_1\).
In case \(m =1\) we have, trivially, \(\alpha_1^M = \alpha_1^{L_1} = e_m \).
When \(m > 1\),
\(u^M = (1/m) e_{m} + [\alpha_{m}]_1^{M \setminus \{m\}}\),
\(u^{L_1} = (1/m) e_{m} + [\alpha_{m}]_1^{L_1 \setminus \{m\}}\)
and \(\alpha_1^M = u^M \, / \, \|u^M\|\),
\(\alpha_1^{L_1} = u^{L_1} \, / \, \|u^{L_1}\|\).

It follows that \(\mathrm{ supp } \, [\alpha_m]_1^{L_1 \setminus \{m\}}\)
is an initial segment of \(M \setminus \{m\}\), and so we infer from the
induction hypothesis applied to \(\alpha_m \), that
\([\alpha_m]_1^{L_1 \setminus \{m\}} = [\alpha_m]_1^{M \setminus \{m\}}\).
Thus \(\alpha_1^M = \alpha_1^{L_1}\) which completes
the case \(n=1\).

Assume now the assertion holds for \(n-1\)
and write \(M = \cup_{i=1}^n
\mathrm{ supp } \,
\alpha_1^{L_i} \cup N \),
where \(\cup_{i=1}^n
\mathrm{ supp } \,
\alpha_1^{L_i} \) is an initial segment of \(M\), which is disjoint
from \(N\).
The induction hypothesis for \(n-1\) yields
\(\alpha_i^M = \alpha_1^{L_i}\) for all \(i < n\).
Hence \(M = \cup_{i=1}^{n-1} \mathrm{ supp } \, \alpha_i^M \cup P\),
where \(P = \mathrm{ supp } \, \alpha_1^{L_n} \cup N\).
It follows from the definition that
\(\alpha_n^M = \alpha_1^P\).
But now the case \(n=1\) guarantees that
\(\alpha_1^P = \alpha_1^{L_n}\)
and the assertion of the lemma is settled.
\end{proof}
{\bf Terminology}. Let \((e_n)\) be a normalized Schauder basic
sequence in a Banach space and let \(\mathcal{F}\) be a regular
family. A finite block basis \(u_1 < \dots < u_m \) of \((e_n)\)
is said to be {\em \(\mathcal{F}\)-admissible} if \(\{\min
\mathrm{ supp } \, u_i :\, i \leq m\} \in \mathcal{F}\). It is
called {\em maximally} \(\mathcal{F}\)-admissible, if
\(\mathcal{F}\) is additionally assumed to be stable and \(\{\min
\mathrm{ supp } \, u_i :\, i \leq m\}\) is a maximal member of
\(\mathcal{F}\).
\begin{Def}
A normalized block basis \((u_n)\) of \((e_n)\)
with \(u_1 < u_2 < ...\) is a \(c_0^\xi\)-spreading model,
if there exists a constant \(C > 0\) such that
\(\|\sum_{i \in F} a_i u_i \| \leq C \max_{i \in F}\, |a_i|\),
for every \(F \in [\mathbb{N}]^{< \infty}\) with
\((u_i)_{i \in F}\) \(S_\xi\)-admissible,
and every choice of scalars \((a_i)_{i \in F}\).
\end{Def}
In what follows we fix a normalized basic sequence
\(\vec{s} = (e_n)\) and a regular and stable family \(\mathcal{F}\).
We abbreviate \(\alpha_n^{\mathcal{F}, \vec{s}, M}\) to
\(\alpha_n^M\).

{\bf Terminology}. Suppose that \(\alpha < \omega_1\) and
\(M \in [\mathbb{N}]\).
An \(\alpha\)-average of \((e_n)\) supported by \(M\),
is any vector of the form \(\alpha_1^L\) for some \(L \in [M]\).

In the sequel we shall make use of the infinite Ramsey theorem
\cite{ell}, \cite{o1} and so we recall its statement.
\([\mathbb{N}]\) is endowed with the
topology of pointwise convergence.
\begin{Thm} \label{ram}
Let \(\mathcal{A}\) be an analytic subset of \([\mathbb{N}]\).
Then there exists \(N \in [\mathbb{N}]\) so that
either \([N] \subset \mathcal{A}\), or
\([N] \cap \mathcal{A} = \emptyset\).
\end{Thm}
Our next result is inspired by an unpublished
result of W.B. Johnson (see \cite{o1}).
\begin{Lem} \label{L13}
Let \(\alpha\) and \(\gamma\) be countable ordinals and suppose
there exists \(N \in [\mathbb{N}]\) such that for every \(M \in
[N]\) there exists a block basis of \(\alpha\)-averages of
\((e_n)\), supported by \(M\), which is a \(c_0^\gamma\)-spreading
model. Then there exist \(M \in [N]\) and a constant \(C > 0\) so
that \(\|\sum_{i=1}^{n_L} \alpha_i^L \| \leq C\), for every \(L
\in [M]\), where \(n_L\) stands for the unique integer satisfying
\(\{\min \mathrm{ supp } \, \alpha_i^L : \, i \leq n_L \}\) is
maximal in \(S_\gamma\).
\end{Lem}
\begin{proof}
Define \(\mathcal{D}_k = \{ L \in [N] : \,
\|\sum_{i=1}^{n_L} \alpha_i^L \| \leq k\}\),
for all \(k \in \mathbb{N}\).
\(\mathcal{D}_k\) is closed in the topology
of pointwise convergence, thanks to Lemma \ref{L12}.
We claim that there exist \(k \in \mathbb{N}\) and
\(M \in [N]\) so that \([M] \subset \mathcal{D}_k\).
The assertion of the lemma clearly follows once this
claim is established.
Were the claim false, then Theorem \ref{ram}
would yield a nested sequence \(M_1 \supset M_2 \supset \dots\)
of infinite subsets of \(N\) such that
\([M_k] \cap \mathcal{D}_k = \emptyset\), for all \(k \in \mathbb{N}\).
Choose an infinite sequence of integers \(m_1 < m_2 < \dots\)
with \(m_i \in M_i\) for all \(i \in \mathbb{N}\).
Set \(M = (m_i)\). Since \(M \in [N]\) our assumptions yield
a block basis \((u_i)\) of \(\alpha\)-averages of \((e_i)\),
supported by \(M\), which is a \(c_0^\gamma\)-spreading model.
Therefore there exists a constant \(C > 0\)
such that \(\|\sum_{i \in F} u_i \| \leq C\),
whenever \((u_i)_{i \in F}\) is \(S_\gamma\)-admissible.
Choose \(k \in \mathbb{N}\) with \(k > C\).
Then choose \(i_0 \in \mathbb{N}\) so that
\(\mathrm{ supp } \, u_i \subset M_k\), for all \(i > i_0\).
If we set \(L = \cup_{i=i_0 + 1}^\infty \mathrm{ supp } \, u_i\),
then \(L \in [M_k]\), and \(\alpha_i^L = u_{i + i_0} \), for
all \( i \in \mathbb{N}\), by Lemma \ref{L12}.
Hence, \(L \notin \mathcal{D}_k\).
However,
\[ \bigl \| \sum_{i=1}^{n_L} \alpha_i^L \bigr \| =
\bigl \| \sum_{i=1}^{n_L} u_{i_0 + i} \bigr \| \leq C < k\]
which is a contradiction.
\end{proof}
\section{Convolution of transfinite averages} \label{S3}
We fix a normalized \(2\)-basic, shrinking sequence
\(\vec{s} = (e_i)\) in some Banach space.
We shall often make use of the following result established
in \cite{o2}: Given \(\epsilon > 0\) there exists
\(M \in [\mathbb{N}]\) such that for every
finitely supported scalar sequence \((a_i)_{i \in M}\)
with \(\|\sum_{i \in M} a_i e_i \| = 1\), we
have \(\max_{i \in M} |a_i| \leq 1 + \epsilon\).
For the rest of this section, we let
\(\mathcal{F} = S_1\). {\em All transfinite
averages of \(\vec{s}\) will be taken with respect to \(\mathcal{F}\)}.
As in the previous section, \(\alpha_n^M\) abbreviates
\(\alpha_n^{\mathcal{F}, \vec{s}, M}\).

The purpose of the present section is to deal with the following problem:
Let \(\alpha\) and \(\beta\) be countable ordinals and suppose that
\((u_i)\) is a block basis of \((\alpha + \beta)\)-averages of \(\vec{s}\).
Does there exist a block basis \((v_i)\) of \(\alpha\)-averages of \(\vec{s}\)
such that \((u_i)\) is a block basis of \(\beta\)-averages of \((v_i)\) ?

It follows directly from the definitions that this is indeed the case
when \(\beta < \omega\). However, if \(\beta\) is an infinite ordinal,
the preceding question has, in general, a negative answer.

In Proposition \ref{Pe1}, we give a partially affirmative answer
to this question which, roughly speaking, states that every
\((\alpha + \beta)\) average of \(\vec{s}\) can be represented
as a finite sum \(\sum_{i=1}^n \lambda_i w_i\), where
\(w_1 < \dots < w_n \) is an \(S_\beta\)-admissible
block basis of \(\alpha\)-averages of \(\vec{s}\) and
\((\lambda_i)_{i=1}^n\) is a sequence of positive scalars
which are almost equal each other. We employ this result in order
to prove the following theorem about transfinite \(c_0\)-spreading
models of \(\vec{s}\), which will in turn be applied in subsequent sections.
In the sequel, when we refer to a block basis
we shall always mean a block basis of \(\vec{s}\). Also all transfinite
averages will be taken with respect to \(\vec{s}\).
\begin{Thm} \label{Te1}
Let \(\alpha\) and \(\beta\) be countable ordinals and \(N \in
[\mathbb{N}]\). Suppose that for every \(P \in [N]\) there exists
\(M \in [P]\) such that no block basis of \(\alpha\)-averages
supported by \(M\) is a \(c_0^\beta\)-spreading model. Then for
every \(P \in [N]\) and \(\epsilon > 0\) there exists \(Q \in
[P]\) with the following property: Every \((\alpha +
\beta)\)-average \(u\) supported by \(Q\) admits a decomposition
\(u = \sum_{i=1}^n \lambda_i u_i\), where \(u_1 < \dots < u_n\) is
a normalized block basis and \((\lambda_i)_{i=1}^n\) is a sequence
of positive scalars such that
\begin{enumerate}
\item There exists \(I \subset \{1, \dots, n\}\) with \((u_i)_{i \in I}\)
\(S_\beta\)-admissible, and such that \(u_i\) is an \(\alpha\)-average for all \(i \in I\),
while \(\|\sum_{i \in \{1, \dots , n\} \setminus I} \lambda_i u_i \|_{\ell_1} < \epsilon\).
\item \(\max_{i \in I} \lambda_i < \epsilon\).
\end{enumerate}
\end{Thm}
Recall that if \(\sum_{i=1}^n a_i e_i\) is a finite linear combination of \(\vec{s}\)
then we denote by \(\|\sum_{i=1}^n a_i e_i\|_{\ell_1} \) the quantity \(\sum_{i=1}^n |a_i|\).
To prove this theorem we shall need to introduce some terminology.
\begin{Def} \label{mD3}
Let \(\alpha\) and \(\beta\) be countable ordinals and \(\epsilon
> 0\). A normalized block \(u\) is said to admit an \((\epsilon,
\alpha, \beta)\)-decomposition, if there exist normalized blocks
\(u_1 < \dots < u_n\) and positive scalars \((\lambda_i)_{i=1}^n\)
with \(u = \sum_{i=1}^n \lambda_i u_i\) and so that the following
conditions are satisfied:
\begin{enumerate}
\item There exists \(I \subset \{1, \dots, n\}\) with \((u_i)_{i \in I}\)
\(S_\beta\)-admissible, and such that \(u_i\) is an \(\alpha\)-average for all \(i \in I\),
while \(\|\sum_{i \in \{1, \dots , n\} \setminus I} \lambda_i u_i \|_{\ell_1} < \epsilon\).
\item \(| \lambda_i - \lambda_j | < \epsilon\) for all \(i\) and \(j\) in \(I\).
\end{enumerate}
\end{Def}
{\bf Terminology}. The quantity \(\max_{i \in I} \lambda_i\) is called the
{\em weight} of the decomposition. If \(u\) is an \((\alpha + \beta)\)-average
admitting an \((\epsilon, \alpha, \beta)\)-decomposition,
\(u = \sum_{i=1}^n \lambda_i u_i\), satisfying \((1)\), \((2)\), above,
and \(I \subset \{1, \dots, n\}\) is as in \((1)\), then every subset of
\(\{\min \, \mathrm{ supp } \, u_i : \, i \in I\}\)
will be called an \((\epsilon, \alpha, \beta)\)-{\em admissible subset of} \(\mathbb{N}\)
{\em resulting from} \(u\).
It is clear that the collection of all \((\epsilon, \alpha, \beta)\)-admissible subsets
of \(\mathbb{N}\) resulting from some (not necessarily the same) \((\alpha + \beta)\)-average
(for some fixed choices of \(\epsilon\), \(\alpha\), \(\beta\)), forms a hereditary family.
\begin{Lem} \label{Le1}
Let \(P \in [\mathbb{N}]\) and \(\epsilon > 0\).
Assume that for every \(L \in [P]\) there exists an \((\alpha + \beta)\)-average
supported by \(L\) which admits an \((\epsilon, \alpha , \beta)\)-decomposition.
Then there exists \(Q \in [P]\) such that every \((\alpha + \beta)\)-average
supported by \(Q\) admits an \((\epsilon, \alpha , \beta)\)-decomposition.
\end{Lem}
\begin{proof}
Let
\[\mathcal{D} = \{ L \in [P]: \, [\alpha + \beta]_1^L \text{ admits an }
(\epsilon, \alpha , \beta)-\text{decomposition} \}.\]
Lemma \ref{L12} yields that \(\mathcal{D}\) is closed in the topology of
pointwise convergence. Theorem \ref{ram} now implies the existence
of some \(Q \in [P]\) such that either \([Q] \subset \mathcal{D}\), or
\([Q] \cap \mathcal{D} = \emptyset\). Our assumptions rule out the second alternative
for \(Q\). Hence \([Q] \subset \mathcal{D}\) which proves the lemma.
\end{proof}
In the next series of lemmas (Lemma \ref{Le3} and Lemma \ref{Le4}), we describe some criteria for embedding a
Schreier family into an appropriate hereditary family of finite subsets
of \(\mathbb{N}\). These criteria, as well as their proofs, are variants of similar results
contained in \cite{amt}, \cite{ag}. We shall therefore
omit the proofs and refer the reader to the aforementioned papers (see for instance
Propositions 2.3.2 and 2.3.6 in \cite{amt}, or Theorems 2.11 and 2.13 in \cite{ag}).
These lemmas will be applied in
the proof of Proposition \ref{Pe1}, which constitutes the main step towards
the proof of Theorem \ref{Te1}.

{\bf Notation}. Let \(\mathcal{F}\) be a family of finite subsets
of \(\mathbb{N}\) and \(M \in [\mathbb{N}]\). Let \(M = (m_i)\)
be the increasing enumeration of \(M\). We set
\(\mathcal{F}(M) = \bigl \{ \{m_i : \, i \in F\}: \, F \in \mathcal{F} \bigr \}\).
Clearly, \(\mathcal{F}(M) \subset \mathcal{F}\) if \(\mathcal{F}\)
is spreading. We also recall that
\(\mathcal{F}[M] = \mathcal{F} \cap [M]^{< \infty}\).
Finally, for every \(L \in [\mathbb{N}]\) and \(\alpha < \omega_1\),
we let \((F_i^\alpha(L))_{i=1}^\infty\) denote the unique decomposition
of \(L\) into successive, maximal members of \(S_\alpha\).
\begin{Lem} \label{Le3}
Suppose that \(1 \leq \xi < \omega_1\),
\(\mathcal{D}\) is a hereditary family of finite subsets of \(\mathbb{N}\)
and \(N \in [\mathbb{N}]\). Assume that for every \(n \in \mathbb{N}\)
and \(P \in [N]\) there exists \(L \in [P]\) such that
\(\cup_{i=1}^n (F_i^\xi(L) \setminus \{\min F_i^\xi(L) \} ) \in \mathcal{D}\).
Then there exists \(M \in [N]\) such that \(S_{\xi + 1} (M) \subset \mathcal{D}\).
\end{Lem}
\begin{Lem} \label{Le4}
Suppose that \(\mathcal{D}\) is a hereditary family of finite
subsets of \(\mathbb{N}\) and \(N \in [\mathbb{N}]\). Let \(\xi <
\omega_1\) be a limit ordinal and let \((\alpha_n)\) be an
increasing sequence of ordinals tending to \(\xi\). Assume there
exists a sequence \(M_1 \supset M_2 \supset \dots\) of infinite
subsets of \(N\) such that \(S_{\alpha_n}(M_n) \subset
\mathcal{D}\), for all \(n \in \mathbb{N}\). Then there exists \(M
\in [N]\) such that \(S_\xi(M) \subset \mathcal{D}\).
\end{Lem}

In the sequel we shall make use of the following permanence property
of Schreier families established in \cite{otw}:
\begin{Lem} \label{l1}
Suppose that \(\alpha < \beta < \omega_1\). Then there exists
\(n \in \mathbb{N}\) such that for every \(F \in S_\alpha\)
with \(n \leq \min F\) we have \(F \in S_\beta\).
\end{Lem}
We shall also make repeated use of the following result from \cite{ano}:
\begin{Lem} \label{lao}
For every \(N \in [\mathbb{N}]\) there exists \(M \in [N]\)
such that for every \(\alpha < \omega_1\) and \(F \in S_\alpha[M]\)
we have \(F \setminus \{\min F\} \in S_\alpha(N)\).
\end{Lem}
Lemma \ref{lao} combined with Proposition 3.2 in \cite{otw} yields the next
\begin{Lem} \label{conv}
Let \(\alpha\) and \(\beta\) be countable ordinals and \(N \in
[\mathbb{N}]\). Then there exists \(M \in [N]\) such that
\begin{enumerate}
\item For every \(F \in S_\beta[S_\alpha][M]\) we have
\(F \setminus \{\min F\} \in S_{\alpha + \beta}\).
\item For every \(F \in S_{\alpha + \beta}[M]\) we
have \(F \setminus \{\min F\} \in S_{\beta}[S_\alpha]\).
\end{enumerate}
\end{Lem}
\begin{Prop} \label{Pe1}
Let \(\alpha\) and \(\beta\) be countable ordinals and \(N \in
[\mathbb{N}]\). Then given \(\epsilon > 0\) and \(P \in [N]\)
there exist \(Q \in [P]\) and \(R \in [Q]\) such that
\begin{enumerate}
\item Every \((\alpha + \beta)\)-average supported by \(Q\) admits
an
\((\epsilon, \alpha , \beta)\) decomposition.
\item For every \(F \in S_\beta[R]\), \(F \setminus \{\min F\}\) is
an \((\epsilon, \alpha , \beta)\)-admissible set resulting from some
\((\alpha + \beta)\)-average supported by \(Q\).
\end{enumerate}
\end{Prop}
\begin{proof}
Fix \(\alpha < \omega_1\). We prove the assertion of the proposition by
transfinite induction on \(\beta\). The case \(\beta = 1\) follows directly
from the definitions since every \((\alpha+1)\)-average admits an
\((\epsilon, \alpha, 1)\)-decomposition. In fact, in this case, we may take
\(Q=P\) and \(R = \{\min \, \mathrm{ supp } \, \alpha_i^P : \, i \in \mathbb{N}\}\)
and check that \((1)\) and \((2)\) hold.

Now let \(\beta > 1\) and suppose the assertion holds for all ordinals smaller than
\(\beta\). Assume first \(\beta\) is a successor ordinal and let \(\beta -1 \)
be its predecessor. Let \(\epsilon > 0\) and \(P \in [N]\) be given and
choose a sequence of positive scalars \((\delta_i)\) such that
\(\sum_i \delta_i < \epsilon / 4\).
Let \(M \in [P]\). The induction hypothesis for \(\beta -1 \) yields infinite
subsets \(R_1 \subset Q_1\) of \(M\) satisfying \((1)\) and \((2)\) for
\((\delta_1, \alpha, \beta -1)\). Choose a maximal member \(F_1 \) of \(S_{\beta -1}\)
with \(F_1 \subset R_1\). We may choose an \((\alpha + \beta -1)\)-average
\(u_1\), supported by \(Q_1\) and such that \(F_1 \setminus \{\min F_1\}\) is
\((\delta_1, \alpha, \beta - 1)\)-admissible resulting from \(u_1\).

Choose \(M_2 \in [M]\) with \(\min M_2 > \max \, \mathrm{ supp } \, u_1\).
Arguing similarly, we choose a maximal member \(F_2\) of \(S_{\beta -1}\)
with \(F_2 \subset M_2\), and an \((\alpha + \beta -1)\)-average
\(u_2\) supported by \(M_2\), which admits a \((\delta_2, \alpha, \beta -1)\)-decomposition
from which \(F_2 \setminus \{\min F_2\}\) is resulting. We continue in this fashion and obtain
a sequence \(F_1 < F_2 < \dots,\) of successive maximal members of \(S_{\beta -1}[M]\),
and a block basis \(u_1 < u_2 < \dots, \) of \((\alpha + \beta -1)\)-averages
supported by \(M\) such that for all \(i \in \mathbb{N}\),
\begin{equation} \label{ep1}
u_i \text{ admits a } (\delta_i, \alpha, \beta -1)-\text{decomposition}.
\end{equation}
\begin{equation} \label{ep2}
F_i \setminus \{\min F_i \} \text{ is } (\delta_i, \alpha, \beta -1)-\text{admissible, resulting from }
u_i.
\end{equation}
We next let, for all \(i \in \mathbb{N}\), \(d_i\) denote the weight of the \((\delta_i, \alpha, \beta -1)\)-decomposition
of \(u_i\), from which \(F_i \setminus \{\min F_i \}\) is resulting.
Clearly, \(d_i \in (0, 3]\).
Therefore, without loss of generality, by passing to a subsequence if necessary,
we may assume that
\begin{equation} \label{ep3}
|d_i - d_j | < \epsilon /4, \text{ for all } i, j \text{ in } \mathbb{N}.
\end{equation}
Now let \(n \in \mathbb{N}\) and choose \(n < i_1 < \dots < i_m \)
such that \((u_{i_k})_{k=1}^m\) is maximally \(S_1\)-admissible.
Set \(u = (\sum_{k=1}^m u_{i_k} ) \, / \, \| \sum_{k=1}^m u_{i_k} \|\).
It is clear that \(u\) is an \((\alpha + \beta)\)-average supported by \(M\).
It is easy to check, using \eqref{ep1} and \eqref{ep3}, that \(u\)
admits an \((\epsilon, \alpha, \beta)\)-decomposition. On the other hand,
\eqref{ep2} implies
that \(\cup_{k=1}^n (F_{i_k} \setminus \{\min F_{i_k}\}) \)
is \((\epsilon, \alpha, \beta)\)-admissible, resulting from \(u\).

Taking in account the stability of \(S_{\beta -1 }\), we conclude the following:
Given \(n \in \mathbb{N}\) and \(M \in [P]\)
\begin{align} \label{ep4}
&\text{There exists an } (\alpha + \beta)-\text{average } u \text{ supported by } M \\
&\text{which admits an } (\epsilon, \alpha, \beta)-\text{decomposition}. \notag
\end{align}
\begin{align} \label{ep5}
&\text{There exists } L \in [M] \text{ such that }
\cup_{i=1}^n (F_i^{\beta -1}(L) \setminus \{ \min F_i^{\beta -1}(L)\}) \\
&\text{ is } (\epsilon, \alpha, \beta)-\text{admissible, resulting from } u. \notag
\end{align}
(Recall that for \(\gamma < \omega_1\),
\((F_i^\gamma(L))_{i=1}^\infty\) denotes the unique decomposition of
\(L\) into consecutive, maximal members of \(S_\gamma\)).

Lemma \ref{Le1} and \eqref{ep4} now yield some \(Q \in [P]\) satisfying \((1)\)
for \((\epsilon, \alpha, \beta)\). Let \(\mathcal{D}\) denote the hereditary
family of the \((\epsilon, \alpha, \beta)\)-admissible subsets of \(Q\)
resulting from some \((\alpha + \beta)\)-average supported by \(Q\). We infer from
\eqref{ep5} that for every \(n \in \mathbb{N}\) and \(M \in [Q]\)
there exists \(L \in [M]\) such that
\(\cup_{i=1}^n (F_i^{\beta -1}(L) \setminus \{ \min F_i^{\beta -1}(L)\})
\in \mathcal{D}\). We deduce from Lemma \ref{Le3} that
\(S_{\beta}(R_0) \subset \mathcal{D}\)
for some \(R_0 \in [Q]\).
Employing Lemma \ref{lao}, we find
\(R \in [R_0]\) such that \(F \setminus \{\min F\} \in \mathcal{D}\),
for all \(F \in S_\beta[R]\). Thus \(Q\) and \(R\) satisfy
\((1)\) and \((2)\) for \((\epsilon, \alpha, \beta)\),
when \(\beta\) is a successor ordinal.

We now consider the case of \(\beta\) being a limit ordinal.
We may choose an increasing sequence of ordinals \((\beta_n)\)
having \(\beta\) as its limit, and such that \((\alpha + \beta_n +1)\)
is the sequence of successor ordinals associated to
the limit ordinal \(\alpha + \beta\).
Let \(\epsilon > 0\) and \(P \in [N]\) be given. Let \(M \in [P]\)
and choose \(m \in M\) with \(1/m < \epsilon /4\).
Next choose \(M_1 \in [M]\) with \(m < \min M_1\)
and such that \(S_{\beta_m}[M_1] \subset S_\beta\) (see Lemma \ref{l1}).
We now apply the induction hypothesis for \(\beta_m\) to obtain
an \((\alpha + \beta_m)\)-average \(v\) supported by \(M_1\) and
admitting an \((\epsilon / 4, \alpha, \beta_m)\)-decomposition.
It is clear that \(u = ((1/m) e_m + v) \, / \, \|(1/m) e_m + v\|\),
is an \((\alpha + \beta)\)-average supported by \(M\) and admitting an
\((\epsilon, \alpha, \beta)\)-decomposition.
Note also that if \(F\) is \((\epsilon / 4, \alpha, \beta_m)\)-admissible
resulting from \(v\), then it is also \((\epsilon, \alpha, \beta)\)-admissible
resulting from \(u\).

It follows now, by lemma \ref{Le1}, that there exists \(Q \in [P]\) such that
\((1)\) holds for \((\epsilon, \alpha, \beta)\).
Next choose positive integers \(k_1 < k_2 < \dots\)
such that \(S_{\beta_n}[k_n, \infty) \subset S_\beta\) (see Lemma \ref{l1}),
for all \(n \in \mathbb{N}\). Successive applications of
the inductive hypothesis applied to each \(\beta_n\) and Lemma \ref{lao},
yield infinite subsets \(Q_1 \supset R_1 \supset Q_2 \supset R_2 \supset \dots\)
of \(Q\) with \(k_n < \min Q_n \) and such that each member of
\(S_{\beta_n}(R_n) \) is an \((\epsilon / 4, \alpha, \beta_n)\)-admissible
set resulting from some \((\alpha + \beta_n)\)-average supported by \(Q_n\),
for all \(n \in \mathbb{N}\).
Let \(\mathcal{D}\) denote the hereditary
family of the \((\epsilon, \alpha, \beta)\)-admissible subsets of \(Q\)
resulting from some \((\alpha + \beta)\)-average supported by \(Q\).
Our preceding argument shows that \(S_{\beta_n} (R_n) \subset \mathcal{D}\),
as long as \(n \in Q\) and \(1/ n < \epsilon /4\).
We deduce now from Lemma \ref{Le4}, that there exists
\(R_0 \in [Q]\) such that \(S_\beta (R_0) \subset \mathcal{D}\).
Once again, Lemma \ref{lao} yields
some \(R \in [R_0]\) with the property
\(F \setminus \{\min F\} \in \mathcal{D}\),
for all \(F \in S_\beta[R]\). Hence, \(Q \supset R\)
satisfy \((1)\) and \((2)\) for \((\epsilon, \alpha, \beta)\),
when \(\beta\) is a limit ordinal. This completes the inductive step
and the proof of the proposition.
\end{proof}
In the proof of Theorem \ref{Te1} we shall need Elton's nearly unconditional
theorem (\cite{e}, \cite{o1}).
\begin{Thm} \label{nun}
Let \((f_i)\) be a normalized weakly null sequence in some Banach space.
There exists a subsequence \((f_{m_i})\) of \((f_i)\) with the following
property: For every \( 0 < \delta \leq 1\) there exists a constant
\(C(\delta) > 0\) such that
\(\|\sum_{i \in F} a_i f_{m_i} \| \leq C(\delta) \|\sum_i a_i f_{m_i}\|\),
for every finitely supported scalar sequence \((a_i)\) in \([-1,1]\)
and every \(F \subset \{i \in \mathbb{N} : \, |a_i| \geq \delta \}\).
\end{Thm}
\begin{proof}[Proof of Theorem \ref{Te1}.]
Let \(P \in [N]\) and \(\epsilon > 0\). Set
\begin{align}
\mathcal{D} = &\{ L \in [P]: [\alpha + \beta]_1^L \text{ admits an } \notag \\
&(\epsilon, \alpha, \beta)-\text{decomposition of weight smaller than } \epsilon \}. \notag
\end{align}
Lemma \ref{L12} yields that \(\mathcal{D}\) is closed in the topology of pointwise
convergence. The theorem asserts that \([Q] \subset \mathcal{D}\), for some \(Q \in [P]\).
Suppose this is not the case and choose, according to Theorem \ref{ram},
\(Q_0 \in [P]\) such that \([Q_0] \cap \mathcal{D} = \emptyset\).
Next choose \(Q_1 \in [Q_0]\) such that no block basis of \(\alpha\)-averages
supported by \(Q_1\) is a \(c_0^\beta\)-spreading model.
Let \(M \in [Q_1]\).
We infer from Proposition \ref{Pe1} that there exist
infinite subsets \(R \subset Q\) of \(M\) such that
\begin{align} \label{ep6}
&\text{Every  } (\alpha + \beta)-\text{average supported by } Q
\text{  admits an } \\
&(\epsilon / 2 , \alpha , \beta)-\text{decomposition}. \notag
\end{align}
\begin{align} \label{ep7}
&\text{If } F \in S_\beta[R],  \text{ then } F \setminus \{\min F\}
\text{ is }
(\epsilon / 2 , \alpha , \beta)-\text{admissible} \\
&\text{resulting from some }
(\alpha + \beta)-\text{average supported by } Q. \notag
\end{align}
Choose a maximal member \(F\) of \(S_\beta[R]\). \eqref{ep6} and
\eqref{ep7} allow us to find normalized blocks \(u_1 < \dots <
u_n\), positive scalars \((\lambda_i)_{i=1}^n\) and \(I \subset
\{1, \dots ,n\}\) such that {\allowdisplaybreaks \begin{align}
\label{ep8}
&\sum_{i=1}^n \lambda_i u_i \text{ is an } (\alpha + \beta)-\text{average supported by } Q, \\
&(u_i)_{i \in I} \text{ is } S_\beta-\text{admissible and }
F \setminus \{ \min F\} \subset \{ \min \, \mathrm{ supp } \, u_i: \, i \in I \}, \notag \\
&u_i \text{ is an } \alpha-\text{average for all } i \in I \text{ and }
\|\sum_{i \in \{1, \dots, n\} \setminus I} \lambda_i u_i \|_{\ell_1} < \epsilon /2, \notag \\
&|\lambda_i - \lambda_j | < \epsilon / 2, \text{ for all } i, j \text{ in } I. \notag
\end{align}}
Since \(\sum_{i=1}^n \lambda_i u_i\) is supported by \(Q \subset Q_0\),
and \([Q_0 ] \cap \mathcal{D} = \emptyset\), we must have that
\(\max_{i \in I} \lambda_i \geq \epsilon\). We deduce from \eqref{ep8}
that
\[\lambda_i \geq \epsilon /2 \, \text{ for all } i \in I.\]
Set \(J_0 = \{ i \in I: \, \min F < \min \, \mathrm{ supp } \, u_i \}\)
and note that \eqref{ep8} implies that
\(F \setminus \{ \min F \} \subset
\{ \min \, \mathrm{ supp } \, u_i: \, i \in J_0 \}\).
It follows now, since \(F\) is maximal in \(S_\beta\), that
\(\{\min F\} \cup \{ \min \, \mathrm{ supp } \, u_i: \, i \in J_0 \} \)
contains a maximal member of \(S_\beta \) as a subset
and therefore, as \(S_\beta\) is stable, there exists an initial segment \(J\) of \(J_0\)
such that
\(\{\min F\} \cup \{ \min \, \mathrm{ supp } \, u_i: \, i \in J \} \)
is a maximal member of \(S_\beta\). Note also that
\(\|\sum_{i \in J} \lambda_i u_i \| \leq 3\).

Summarizing, given \(M \in [Q_1]\) we found a block basis of
\(\alpha\)-averages \(v_1 < \dots < v_k\), supported by \(M\),
\(m \in M\) with \(m < \min \, \mathrm{ supp } \, v_1\), and scalars
\((\mu_i)_{i=1}^k \) in \([\epsilon / 2, 2]\) so that
\begin{equation} \label{ep9}
\{m\} \cup \{ \min \, \mathrm{ supp } \, v_i: \, i \leq k\} \text{ is maximal in }
S_\beta \text{ and }
\|\sum_{i=1}^k \mu_i v_i \| \leq 3.
\end{equation}
Define
\begin{align}
\mathcal{D}_1 = &\bigl \{ L \in [Q_1]: \, \exists \, (\mu_i)_{i=1}^k \subset [\epsilon / 2 , 2], \,
\|\sum_{i=1}^k \mu_i \alpha_i^{L \setminus \{ \min L \}}\| \leq 3, \text{ and } \notag \\
&\{\min L\} \cup \{ \min \, \mathrm{ supp } \, \alpha_i^{L \setminus \{ \min L \}} : \,
i \leq k \} \text{ is maximal in } S_\beta \bigr \}. \notag
\end{align}
Lemma \ref{L12} and the stability of \(S_\beta\) yield that \(\mathcal{D}_1\)
is closed in the topology of pointwise convergence.
We now infer from \eqref{ep9} that every \(M \in [Q_1]\) contains some
\(L \in \mathcal{D}_1\) as a subset. Thus, we deduce from Theorem \ref{ram}
that there exists \(M_0 \in [Q_1]\) with \([M_0] \subset \mathcal{D}_1\).

Now let \(L \in [M_0]\) and denote by \(n_L\) the unique integer such that
\((\alpha_i^L)_{i=1}^{n_L}\) is maximally \(S_\beta\)-admissible.
Because \(L \in \mathcal{D}_1\), we must have that
\begin{align} \label{ep10}
&\|\sum_{i=1}^{n_L} \mu_i \alpha_i^L \| \leq 4,
\text{ for some choice of scalars } \\
&(\mu_i)_{i=1}^{n_L} \text{ in
the interval} \, [\epsilon / 2, 2]. \notag
\end{align}
Set \(g_i = \alpha_i^{M_0}\), for all \(i \in \mathbb{N}\).
Then \((g_i)\) is a normalized weakly null sequence, as \(\vec{s}\)
is assumed to be shrinking.
Theorem \ref{nun}
now yields a constant \(C > 0\) and a subsequence of \((g_i)\)
(which, for clarity, is still denoted by \((g_i)\)), such that
\[ \| \sum_{i \in G} a_i g_i \| \leq C \| \sum_{i=1}^\infty a_i g_i \|,\]
for every finitely supported scalar sequence \((a_i)\) in \([-2,2]\)
and every \(G \subset \{i \in \mathbb{N}: \, |a_i| \geq \epsilon /2\}\).
It follows from this, Lemma \ref{L12} and \eqref{ep10} that, whenever \(F \in [\mathbb{N}]^{< \infty}\)
is so that \((g_i)_{i \in F}\) is maximally \(S_\beta\)-admissible, then
we have some choice of scalars \((\mu_i)_{i \in F}\) in \([\epsilon / 2, 2]\)
such that
\[\| \sum_{i \in F} \sigma_i \mu_i g_i \| \leq 8C ,\]
for every choice of signs \((\sigma_i)_{i \in F}\).
We conclude from the above, that some subsequence of \((g_i)\)
is a \(c_0^\beta\)-spreading model. Lemma \ref{L12}
finally implies that there is some \(L \in [M_0]\) (and thus \(L \in [Q_1]\))
such that \((\alpha_i^L)\) is a \(c_0^\beta\)-spreading model,
contradicting the choice of \(Q_1\).
Therefore, we must have that
\([Q] \subset \mathcal{D}\), for some \(Q \in [P]\),
and the proof of the theorem is now complete.
\end{proof}
\section{Transfinite averages of weakly null sequences
in \(C(K)\) equivalent to the unit vector basis of \(c_0\)} \label{S4}
In this section we present the following
\begin{Thm} \label{mT}
Let \(K\) be a compact metric space and let \((f_n)\) be a
normalized, basic sequence in \(C(K)\). Suppose that there exist
\(M \in [\mathbb{N}]\) and a summable sequence of positive scalars
\((\epsilon_n)\) such that for all \(t \in K\), the set \(\{n \in
M : \, |f_n(t)| \geq \epsilon_n \}\) is finite. Then there exist
\(\xi < \omega_1\) and a block basis of \(\xi\)-averages of
\((f_n)\) equivalent to the unit vector basis of \(c_0\).
\end{Thm}
(Note that all transfinite averages of \((f_n)\) are considered
with respect to \(\mathcal{F} = S_1\).)
\begin{remark}
The hypotheses in Theorem \ref{mT} imply that
\(\sum_{n \in M} |f_n(t)|\) is a convergent series,
for all \(t \in K\). It follows then from Rainwater's theorem \cite{rw},
that every normalized block basis of \((f_n)_{n \in M}\)
is weakly null and therefore, the
subsequence \((f_n)_{n \in M}\) of \((f_n)\) is
shrinking. Moreover, the convergence of the series
\(\sum_{n \in M} |f_n(t)|\)
for all \(t \in K\), implies that some block basis
of \((f_n)_{n \in M}\) is equivalent to the unit vector
basis of \(c_0\). This is a special case of a famous
result, due to J. Elton \cite{e2}, which states that
if \((x_n)\) is a normalized basic sequence in some
Banach space and the series \(\sum_n |x^*(x_n)|\)
converges for every extreme point \(x^*\) in the
ball of \(X^*\), then some block basis of \((x_n)\)
is equivalent to the unit vector basis of \(c_0\).
An alternate proof of this special case of Elton's
theorem is given in \cite{hor}. See also
\cite{fo}, \cite{an} for related results.
We wish to indicate however, as our next corollary shows, that this special case
of Elton's theorem is also a consequence of Theorem
\ref{mT}. Hence, our result may be viewed as a quantitative version
of this special case of Elton's theorem.
\end{remark}
\begin{Cor} \label{elton}
Let \((f_n)\) be a normalized basic sequence
in \(C(K)\) such that \(\sum_n |f_n(t)|\) is a convergent series,
for all \(t \in K\). Then there exist \(\xi < \omega_1\)
and a block basis of \(\xi\)-averages of \((f_n)\)
equivalent to the unit vector basis of \(c_0\).
\end{Cor}
The proof is given at the end of this section.

The ordinal \(\xi\) that appears in the conclusion of Theorem \ref{mT},
is related to the complexity of the compact family
\(\{F \in [M]^{< \infty} : \, \exists \, t \in K \text{ with }
|f_n(t)| \geq \epsilon_n, \, \forall \, n \in F\}\).
It follows from Corollary \ref{C0}, that every normalized weakly null sequence
in \(C(K)\), for \(K\) a countable compact metric space, admits a subsequence satisfying the hypotheses
of Theorem \ref{mT}. Moreover, if \(K\) is homeomorphic to
\([1, \omega^{\omega^\alpha}]\), for some \(\alpha < \omega_1\), then as is shown
in Corollary \ref{mC},
the ordinal \(\xi\) in the conclusion of Theorem \ref{mT} can be taken
not to exceed \(\alpha\).

We shall next describe how to obtain the ``optimal'' \(\xi\) satisfying the conclusion
of Theorem \ref{mT}.

{\bf The following conventions hold throughout this section}.
\(K\) is a compact metric space and \(\vec{s} = (f_n)\)
is a normalized shrinking basic sequence
in \(C(K)\). We shall assume, without loss
of generality, by passing to a subsequence if necessary,
that \(\vec{s}\) is \(2\)-basic.
We let
\(\mathcal{F} = S_1\). {\em All transfinite
averages of \(\vec{s}\) will be taken with respect to \(\mathcal{F}\)}.
As in the previous section, \(\alpha_n^M\) abbreviates
\(\alpha_n^{\mathcal{F}, \vec{s}, M}\).
In the sequel, when we refer to a block basis
we shall always mean a block basis of \(\vec{s} = (f_n)\). Also all transfinite
averages will be taken with respect to \(\vec{s}\).
\begin{Def} \label{mD1}
\begin{enumerate}
\item Given \(N \in [\mathbb{N}] \) and \( 1 \leq \alpha < \omega_1\),
we say that \(N\) is \(\alpha\)-large, if for every \(\beta < \alpha\)
and \(M \in [N]\) there exists \(L \in [M]\) such that no block basis
of \(\beta\)-averages supported by \(L\) is a \(c_0^\gamma\)-spreading model,
where \(\beta + \gamma = \alpha\).
\item Given \(N \in [\mathbb{N}] \) set
\(\xi^N = \sup \{ \alpha < \omega_1 : \, \exists \text{ an } \alpha-\text{large} \,
M \in [N] \}\).
Put \(\xi^N = 0\), if this set is empty. Finally put
\(\xi^0 = \min \{ \xi^N : \, N \in [\mathbb{N}]\}\).
\end{enumerate}
\end{Def}
Note that if \(\xi^0 = \xi^{N_0}\) for some \(N_0 \in [\mathbb{N}]\),
then \(\xi^L = \xi^0\), for all \(L \in [N_0]\). In fact,
if \(1 \leq \xi^0 < \omega_1\),
then every infinite subset of \(N_0\) is \(\xi^0\)-large.
\begin{Prop} \label{mP}
Suppose that \(\xi^N < \omega_1\), for some \(N \in [\mathbb{N}]\).
Then there exists a block basis of \(\xi^N\)-averages, supported by
\(N\), which is equivalent to the unit vector basis of \(c_0\).
\end{Prop}
We postpone the proof and observe that if \(\xi^0 < \omega_1\)
and \(\xi^0 = \xi^{N_0}\), then Proposition \ref{mP} yields that every
infinite subset of \(N_0\) supports a block basis of \(\xi^0\)-averages,
equivalent to the unit vector basis of \(c_0\) and, moreover, it follows
by our preceding comments, that
\(\xi^0\) is the smallest ordinal with this property.
Therefore, the optimality of \(\xi^0\) is considered in this sense.
In order to prove Theorem \ref{mT}, we need to introduce some more
notation and terminology.
\begin{Def} \label{mD2}
\begin{enumerate}
\item Let \(\beta < \alpha < \omega_1\),
\(p \in \mathbb{N}\) and \(\epsilon > 0\).
An \(\alpha\)-average \(u = \sum_i a_i f_i\),
is said to be \((\beta, p, \epsilon)\)-large, if
for every choice \(I_1 < \dots < I_k\)
of \(k\) consecutive members of \(S_\beta\), \(k \leq p\), and all
\(t \in K\), we have
\(|\sum_{i \in I} a_i f_i(t) | \leq  \epsilon
+ \sum_{i \notin I } a_i |f_i(t) |\),
where \(I = \cup_{j=1}^k I_j\).
\item Let \(N \in [\mathbb{N}] \), \( 1 \leq \alpha < \omega_1\)
We say that \(N\) is \(\alpha\)-nice
if for every \(\beta < \alpha\), every \(M \in [N]\),
every \(p \in \mathbb{N}\) and all \(\epsilon > 0\),
there exists an \(\alpha\)-average
supported by \(M\) which is \((\beta, p, \epsilon)\)-large.
\end{enumerate}
\end{Def}
The main step for proving Theorem \ref{mT} is
\begin{Thm} \label{mT2}
Suppose that \(N \in [\mathbb{N}]\)
is \(\alpha\)-large, for some \( 1 \leq \alpha < \omega_1\).
Then \(N\) is \(\alpha\)-nice.
\end{Thm}
We postpone the proof in order to give the
\begin{proof}[Proof of Theorem \ref{mT}]
Let
\[\mathcal{G}= \{F \in [\mathbb{N}]^{< \infty} : \, \exists \, t \in K \text{ with }
|f_n(t)| \geq \epsilon_n, \, \forall \, n \in F\}. \]
Clearly, \(\mathcal{G}\) is hereditary. The compactness of \(K\) and
our assumptions, imply that \(\mathcal{G}[M]\) is compact in the topology
of pointwise convergence.
It follows that there is a countable
ordinal \(\zeta\) such that \(\mathcal{G}[M]^{(\zeta)}\)
is finite.
Write \(\zeta = \omega^\gamma k + \eta\), for some
\(k \in \mathbb{N}\) and \(\eta < \omega^\gamma\).
We infer now by the result of \cite{g}, that
there exists \(N \in [M]\) with the property
\(\mathcal{G}[N] \subset S_{\gamma + 1}\).

We claim that \(\xi^N \leq \gamma + 1\) (see Definition \ref{mD1}).
Indeed, were this claim false, we would choose \(P \in [N]\)
and a countable ordinal \(\beta > \gamma + 1\) such that
\(P\) is \(\beta\)-large. Theorem \ref{mT2} then
yields \(P\) is \(\beta\)-nice
(see Definition \ref{mD2}).
Next, let \(\epsilon > 0\) and choose \(Q \in [P]\) such that
\(\sum_{n \in Q} \epsilon_n < \epsilon / 12\).
Since \(\gamma + 1 < \beta\) and \(P\) is
\(\beta\)-nice,
there exists a \(\beta\)-average \(u = \sum_i a_i f_i\),
supported by \(Q\) which is \((\gamma +1 , 1, \epsilon/2)\)-large.
This means
\[\bigl |\sum_{i \in I} a_i f_i(t) \bigr | \leq \epsilon /2
+ \sum_{ i \notin I} a_i |f_i(t)|, \]
for all \(t \in K\)
and every \(I \in S_{\gamma + 1}\).
Given \(t \in K\), put
\(\Lambda_t = \{ n \in \mathbb{N} : \, |f_n(t)| \geq \epsilon_n \}\).
Note that \(u\) is supported by \(N\) and so
\(\Lambda_t \cap \mathrm{ supp } \, u \in S_{\gamma + 1}\),
for all \(t \in K\), as
\(\Lambda_t \cap \mathrm{ supp } \, u \in \mathcal{G}[N]\).
Taking in account that \(\|u\| = 1\), we have \(0 \leq a_i \leq 3\),
for all \(i \in \mathbb{N}\). Hence,
\begin{align}
|u(t)| &\leq \biggl | \sum_{ i \in \Lambda_t \cap \mathrm{ supp } \, u} a_i f_i(t) \biggr |
+ \biggl |\sum_{ i \notin \Lambda_t } a_i f_i(t) \biggr | \notag \\
&\leq \epsilon /2 + 2 \sum_{ i \notin \Lambda_t} a_i |f_i(t)| \notag \\
&< \epsilon /2 + 6 \epsilon /12 = \epsilon, \notag
\end{align}
for all \(t \in K\). Since \(\epsilon\) was arbitrary, we have
reached a contradiction. Therefore, our claim holds. In
particular, \(\xi^N < \omega_1\) and the assertion of the theorem
is a consequence of Proposition \ref{mP}.
\end{proof}
\begin{Cor} \label{mC}
Let \((f_n)\) be a normalized weakly null sequence in
\(C(\omega^{\omega^\xi})\), \(\xi < \omega_1\). Then there exist \(\alpha \leq \xi\)
and a block basis of \(\alpha\)-averages of \((f_n)\) equivalent
to the unit vector basis of \(c_0\).
\end{Cor}
\begin{proof}
Set \(K = [1,\omega^{\omega^\xi}]\).
Corollary \ref{C0} yields \(M \in [\mathbb{N}]\) and
a summable sequence of positive scalars \((\epsilon_n)\) such that
for all \(t \in K\) the set
\(\{n \in M : \, |f_n(t)| \geq \epsilon_n \}\)
belongs to \(S_\xi^{+}\). In particular,
\(\Lambda_t \cap M\) is the union of two consecutive members of
\(S_\xi\). The argument in the proof of Theorem \ref{mT}
shows that \(\xi^M \leq \xi\). The assertion of the corollary
now follows from Proposition \ref{mP}.
\end{proof}
We shall now give the proof of Proposition \ref{mP}.
This requires two lemmas.
\begin{Lem} \label{Le5}
Suppose that \(1 \leq \alpha < \omega_1\).
Let \(m < n\) in \(\mathbb{N}\) and \(F \in [\mathbb{N}]^{< \infty}\)
with \(n < \min F\) be such that \(\{n\} \cup F \) is a maximal member
of \(S_\alpha\). Then \(\{m\} \cup F \notin S_\alpha\).
\end{Lem}
\begin{proof}
We use transfinite induction on \(\alpha\). When \(\alpha =1\),
we must have that \(|F| = n-1\), in order for \(\{n\} \cup F \)
be maximal in \(S_1\). Hence,
\(|\{m\} \cup F | = n > m = \min (\{m\} \cup F )\).
Thus the assertion of the lemma holds in this case.

Next assume the assertion holds for all ordinals smaller than
\(\alpha\) (\(\alpha > 1\)).
Suppose first \(\alpha\) is a limit ordinal and let \((\alpha_n)\)
be the sequence of successor ordinals associated to \(\alpha\).
Since \(\{n\} \cup F \) is maximal in \(S_\alpha\), we have that
\(\{n\} \cup F \) is maximal in \(S_{\alpha_k}\),
for all \(k \leq n\) such that
\(\{n\} \cup F \in S_{\alpha_k}\).
Suppose we had \(\{m\} \cup F \in S_\alpha\).
Then there is some \(k \leq m\) such that
\(\{m\} \cup F \in S_{\alpha_k}\).
We infer from the spreading property of
\(S_{\alpha_k}\), as \(m < n\), that
\(\{n\} \cup F \in S_{\alpha_k}\).
Therefore, \(\{n\} \cup F \) is maximal in \(S_{\alpha_k}\).
The induction hypothesis applied on \(\alpha_k\)
now yields \(\{m\} \cup F \notin S_{\alpha_k}\),
a contradiction which proves the assertion when
\(\alpha\) is a limit ordinal.

We now assume \(\alpha\) is a successor ordinal, say
\(\alpha = \beta + 1\).
Since \(\{n\} \cup F \) is maximal in \(S_\alpha\),
there exist \(F_1 < \dots < F_n\),
successive maximal members of \(S_\beta\)
such that \(\{n\} \cup F = \cup_{i=1}^n F_i\)
(see \cite{g}).
We shall assume \(m > 1\) or else the assertion holds
since \(\{1\}\) is maximal in every Schreier family
and \(F \ne \emptyset\).
Note that the induction hypothesis on \(\beta\)
implies that
\(G_1 = \{m\} \cup (F_1 \setminus \{n\}) \notin S_\beta\).
It follows, as \(S_\beta\) is stable, that
\(G_1\) contains a maximal member \(H_1\) of
\(S_\beta\) as an initial segment,
and so we may write \(G_1 = H_1 \cup H_2\)
with \(H_2 \ne \emptyset\). Of course, \(m = \min H_1\).
Set \(H = H_1 \cup \cup_{i=2}^m F_i\).
Then \(H\) is maximal in \(S_\alpha\).
This completes the proof of the lemma since
\(H\) is a proper subset of
\(\{m\} \cup F\).
\end{proof}
\begin{Lem} \label{Le6}
Let \(P \in [\mathbb{N}]\), \(\beta \leq \alpha < \omega_1\)
and \(\tau < \omega_1\). Assume that every block basis
of \(\beta\)-averages supported by \(P\) is a \(c_0^\gamma\)-spreading model,
where \(\beta + \gamma = \alpha\), while every block basis
of \(\alpha\)-averages supported by \(P\) is a \(c_0^\tau\)-spreading model.
Then there exists \(Q \in [P]\) such that every block basis of
\(\beta\)-averages supported by \(Q\) is a \(c_0^{\gamma + \tau}\)-spreading model.
\end{Lem}
\begin{proof}
We assume that both \(\gamma\) and \(\tau\) are greater than or equal
to \(1\), or else the assertion of the lemma is trivial.
We also assume, without loss of generality thanks to Lemma \ref{L13},
that there exists a constant \(C > 0\) such that every block basis
of \(\beta\)-averages (resp. \(\alpha\)-averages)
supported by \(P\) is a \(c_0^\gamma\)
(resp. \(c_0^\tau\))-spreading model
with constant \(C\).
We shall further assume, without loss of generality thanks
to Lemma \ref{conv},
that for every
\(F \in S_{\gamma + \tau}[P]\)
we have \(F \setminus \{\min F\} \in S_\tau[S_\gamma]\).

Let \(M \in [P]\). Choose a sequence of positive scalars
\((\delta_i)\) with \(\sum_i \delta_i < 1/(4C)\).
We apply Proposition \ref{Pe1}, successively, to obtain
the following objects:
\begin{enumerate}
\item A maximally \(S_\tau\)-admissible block basis
\(v_1 < \dots < v_n\) of \(\alpha\)-averages, supported by \(M\),
with \(\min M < \min \mathrm{ supp } \, v_1\).
\item Successive, maximal members \(F_1 < \dots < F_n\)
of \(S_\gamma[M]\) such that \(\max \mathrm{ supp } \, v_i < \min F_{i+1}\),
for all \(i < n\).
\item Successive finite subsets of \(\mathbb{N}\)
\(J_1 < \dots < J_n\) such that for each \(i \leq n\),
there exist a normalized block basis \((u_j)_{j \in J_i}\),
a subset \(I_i \) of \(J_i\) and positive scalars \((\lambda_j)_{j \in J_i}\)
which satisfy the following properties:
\end{enumerate}
\begin{equation} \label{meq1}
v_i = \sum_{j \in J_i} \lambda_j u_j, \text{ and }
\biggl \| \sum_{j \in J_i \setminus I_i} \lambda_j u_j \biggr \|_{\ell_1} < \delta_i.
\end{equation}
\begin{align} \label{meq2}
&(u_j)_{j \in I_i} \text{ is an } S_\gamma-\text{admissible block basis of }
\beta-\text{averages} \\
&\text{and } |\lambda_r - \lambda_s | < \delta_i, \text{ for all }
r, s \text{ in } I_i. \notag
\end{align}
\begin{equation} \label{meq3}
F_i \setminus \{\min F_i \} \subset \{\min \mathrm{ supp } \, u_j : \, j \in I_i\}.
\end{equation}
Our assumptions yield that \(\|\sum_{i=1}^n v_i \| \leq C\) and that
\[ 1 - \delta_i \leq \biggl \| \sum_{j \in I_i} \lambda_j u_j \biggr \| \leq
C \max_{j \in I_i} \lambda_j, \text{ for all } i \leq n.\]
\eqref{meq2} now implies
\begin{equation} \label{meq4}
1/(2C) \leq \lambda_j \leq 3, \text{ for all } j \in I_i \text{ and } i \leq n.
\end{equation}
We also obtain from \eqref{meq1} that
\begin{equation} \label{meq5}
\biggl \| \sum_{i=1}^n \sum_{j \in I_i} \lambda_j u_j \biggr \|
\leq C + \sum_{i=1}^n \delta_i < 2C.
\end{equation}
We next observe that for all \(i < n\) and \(j_0 \in I_i\),
\(\{\min \mathrm{ supp } \, u_{j_0} \} \cup
\{\min \mathrm{ supp } \, u_j : \, j \in I_{i+1} \}
\notin S_\gamma\).
This is so since
\(F_{i+1} \setminus \{\min F_{i+1} \} \subset \{\min \mathrm{ supp } \, u_j : \, j \in I_{i+1}\}\),
(by \eqref{meq3}), \( \max \mathrm{ supp } \, v_i < \min F_{i+1}\),
and thus, as a consequence of Lemma \ref{Le5}, we have that
\(\{\min \mathrm{ supp } \, u_{j_0} \} \cup (F_{i+1} \setminus \{\min F_{i+1} \} )
\notin S_\gamma\).

It follows from this that for all \(i \leq n\) there exists an initial
segment \(I_i^{*}\) of \(I_i\) (possibly, \(I_i^{*} = \emptyset\))
with \(\max I_i^{*} < \max I_i\), such that
\(\{\min \mathrm{ supp } \, u_j : \, j \in (I_i \setminus I_i^{*}) \cup I_{i+1}^{*}\}\)
is a maximal member of \(S_\gamma\), for all \(i < n\).
Note that \(I_1^{*} = \emptyset\).

Set \(T_i = (I_i \setminus I_i^{*}) \cup I_{i+1}^{*}\), for
all \(i < n\). Then \((u_j)_{j \in T_i}\)
is maximally \(S_\gamma\)-admissible for all \(i < n\).
We also infer from \eqref{meq4} and \eqref{meq5} that
\[ \biggl \| \sum_{j \in \cup_{i < n} T_i} \lambda_j u_j \biggr \| \leq 4C, \,
\lambda_j \in [1/(2C), 3], \text{ for all }
j \in \cup_{i < n} T_i.\]
Note also that \(\min \mathrm{ supp } \, u_{\min T_i} <
\min \mathrm{ supp } \, v_{i+1}\), for all \(i < n\).
Since \(\min M < \min \mathrm{ supp } \, v_1\) and
\((v_i)_{i=1}^n \) is maximally \(S_\tau\)-admissible,
Lemma \ref{Le5} and the spreading property of
\(S_\tau\), yield that
\(\{ \min M\} \cup
\{ \min \mathrm{ supp } \, u_{\min T_i} : \, i < n\}\)
is not a member of \(S_\tau\).
Hence, by the stability of \(S_\tau\), there exists \(m < n\) such that
\(\{ \min M\} \cup
\{ \min \mathrm{ supp } \, u_{\min T_i} : \, i \leq m \}\)
is a maximal member of \(S_\tau\).
Note also that
\(\| \sum_{j \in \cup_{i \leq m} T_i} \lambda_j u_j \| \leq 4C\)
and
\(\lambda_j \in [1/(2C), 3]\), for all
\(j \in \cup_{i \leq m} T_i \).

Summarizing, given \(M \in [P]\) there exists a maximally
\(S_\tau[S_\gamma]\)-admissible block basis \((u_i)_{i=1}^k\)
of \(\beta\)-averages, supported by \(M\), and
scalars \((\lambda_i)_{i=1}^k\) in \([1/(2C), 3]\)
such that \(\|\sum_{i=1}^k \lambda_i u_i \| \leq 5C\).
Given \(L \in [P]\) let \(n_L\) denote the unique integer such that
\((\beta_i^L)_{i=1}^{n_L}\) is maximally \(S_\tau[S_\gamma]\)-admissible.
Define
\[
\mathcal{D} = \biggl \{ L \in [P]: \, \exists \, (\lambda_i)_{i=1}^{n_L} \subset [1/(2C) , 3], \,
\bigl \|\sum_{i=1}^{n_L} \lambda_i \beta_i^L \bigr \| \leq 5C \biggr \}.\]
Lemma \ref{L12} and the stability of \(S_\tau[S_\gamma]\) yield that \(\mathcal{D}\)
is closed in the topology of pointwise convergence.
We infer from our preceding discussion, that every \(M \in [P]\) contains some
\(L \in \mathcal{D}\) as a subset. Thus, we deduce from Theorem \ref{ram}
that there exists \(M_0 \in [P]\) with \([M_0] \subset \mathcal{D}\).
Arguing as in the last part of the proof of Theorem \ref{Te1},
using Theorem \ref{nun}
and our assumptions on \(P\),
we obtain a block basis of \(\beta\)-averages which is a
\(c_0^{\gamma + \tau}\)-spreading model.
The assertion of the lemma now follows from Lemma \ref{L13}.
\end{proof}
\begin{proof}[Proof of Proposition \ref{mP}]
To simplify our notation, let us write \(\xi\) instead of \(\xi^N\).
We assert that for every \(M \in [N]\) and all
\(\beta < \omega_1\) there exists a block basis of \(\xi\)-averages
supported by \(M\) which is a \(c_0^\beta\)-spreading model.
Once this is accomplished, the proposition will
follow from the Kunen-Martin boundedness principle (see \cite{d}, \cite{ke}).
To see this, let \(N \in [\mathbb{N}]\). Given \(n \in \mathbb{N}\),
let \(\mathcal{T}_n^N \) denote the family of those finite subsets of \(N\)
that are initial segments of sets of the form
\(\cup_{i=1}^k \mathrm{ supp } \, \xi_i^L\), for some \(k \in \mathbb{N}\)
and \(L \in [N]\)
such that \(\|\sum_{i=1}^k \xi_i^L \| \leq n\).
We claim there is some \(n \in \mathbb{N}\)
so that \(\mathcal{T}_n^N\) is not compact in the topology
of pointwise convergence.
Otherwise, the Mazurkiewicz-Sierpinski theorem \cite{ms},
yields \(\zeta < \omega_1\) so that
\(\mathcal{T}_n^N\) is homeomorphic to a subset of
\([1, \omega^{\omega^\zeta}]\), for all \(n \in \mathbb{N}\).
We may now choose, according to our assertion combined with Lemma \ref{L13},
some \(L_0 \in [N]\) and \(n \in \mathbb{N}\) such that
\((\xi_i^L)_{i=1}^\infty\) is a \(c_0^{\zeta +1}\)-spreading model
with constant \(n\), for all \(L \in [L_0]\).
It follows from this that for all \(L \in [L_0]\),
\(\cup_{i=1}^{n_L} \mathrm{ supp } \, \xi_i^L \in \mathcal{T}_n^N\),
where \(n_L\) stands for the unique integer such that
\((\xi_i^L)_{i=1}^{n_l}\) is maximally \(S_{\zeta +1}\)-admissible.
Since \(S_\alpha\) is homeomorphic to \([1, \omega^{\omega^\alpha}]\) for
all \(\alpha < \omega_1\) (see \cite{aa}), this implies that
\(S_{\zeta + 1}\) is homeomorphic to a subset of
\([1, \omega^{\omega^\zeta}]\) which is absurd.
Hence, indeed, there is some \(n \in \mathbb{N}\) with
\(\mathcal{T}_n^N\) non-compact. Subsequently, there exists
\(M \in [N]\), \(M = (m_i)\), such that
\(\{m_1, \dots, m_k \} \in \mathcal{T}_n^N\), for
all \(k \in \mathbb{N}\). We now infer from Lemma \ref{L12},
that \(\|\sum_{i=1}^k \xi_i^M \| \leq n\), for all \(k \in \mathbb{N}\).
Using an argument based on Theorem \ref{ram}, similar to that in the proof of Lemma
\ref{L13}, we conclude that some block basis of \(\xi\)-averages is equivalent to the
unit vector basis of \(c_0\).

We shall next prove our initial assertion by transfinite induction on \(\beta\).
The assertion is trivial for \(\beta = 0\). Assume \(\beta \geq 1\)
and that the assertion holds for all \(M\in [N]\) and all ordinals
smaller than \(\beta\) yet, for some \(P \in [N] \)
there is no block basis of \(\xi\)-averages, supported by \(P\),
which is a \(c_0^\beta\)-spreading model.
We now show that \(P\) is \((\xi + \beta)\)-large which, of course,
is absurd.

To see this, first consider an ordinal
\(\gamma < \xi\) and let \(M \in [P]\). Write \(\xi = \gamma + \delta\).
We claim that there exists \(L \in [M]\) such that no block basis
of \(\gamma\)-averages supported by \(L\) is a
\(c_0^{\delta + \beta}\)-spreading model (note that
\(\gamma + (\delta + \beta) = \xi + \beta\)).
Were this claim false, then Lemma \ref{L13} would yield a constant
\(C > 0\) and \(L_0 \in [M]\) such that, every block basis of
\(\gamma\)-averages supported by \(L_0\) is a
\(c_0^{\delta + \beta}\)-spreading model
with constant \(C\). By employing Lemma \ref{conv}
we may assume, without loss of generality, that for all
\(F \in S_\beta[S_\delta]\), \(F \subset L_0\),
we have \(F \setminus \{ \min F \} \in S_{\delta + \beta}\).
But now, we shall exhibit a block basis of \(\xi\)-averages
supported by \(L_0\) (and thus also by \(P\)), which is
a \(c_0^\beta\)-spreading model. Indeed, as \(\xi = \gamma + \delta\),
we may apply Proposition \ref{Pe1}, successively, to obtain block bases
\(u_1 < u_2 < \dots \) and \(v_1 < v_2 < \dots \)
consisting of \(\xi\) and \(\gamma\)-averages, respectively,
both supported by \(L_0\);
A sequence of positive scalars \((\lambda_i)\) and a sequence
\(F_1 < F_2 < \dots\) of successive finite subsets of \(\mathbb{N}\)
so that the following requirements are satisfied:
\begin{enumerate}
\item \(\|u_i - \sum_{ j \in F_i} \lambda_j v_j \| < \epsilon_i\),
for all \(i \in \mathbb{N}\).
\item \((v_j)_{j \in F_i}\) is \(S_\delta\)-admissible
and \(\mathrm{ supp } \, v_j \subset \mathrm{ supp } \, u_i \),
for all \(j \in F_i\) and \(i \in \mathbb{N}\).
\end{enumerate}
In the above, \((\epsilon_i)\) is a summable sequence of positive scalars.
Since \((\lambda_j)_{j \in \cup_i F_i}\) is bounded and
\((v_i)\) is a \(c_0^{\delta + \beta}\)-spreading model, our assumptions on \(L_0\)
readily imply that \((u_i)\) is a block basis of \(\xi\)-averages supported by \(P\)
which is a \(c_0^\beta\)-spreading model. This contradicts the choice of \(P\).
Therefore our claim holds.

Next, let \(M \in [P]\), \(\gamma < \beta\) and write \( \beta = \gamma + \delta\).
Note that \(\xi + \beta = (\xi + \gamma ) + \delta\). We now claim that
there exists \(L \in [M]\) such that no block basis of \((\xi + \gamma)\)-averages
supported by \(L\) is a \(c_0^\delta\)-spreading model.
If that were not the case then, thanks to Lemma \ref{L13},
there would exist \(L_0 \in [M]\) such that every block basis of
\((\xi + \gamma)\)-averages
supported by \(L_0\) is a \(c_0^\delta\)-spreading model.

Since \(\gamma < \beta\), the induction hypothesis combined with Lemma
\ref{L13} implies the existence of some \(L_1 \in [L_0]\)
such that every block basis of \(\xi\)-averages supported by \(L_1\)
is a \(c_0^\gamma\)-spreading model. We deduce from Lemma \ref{Le6}
that some block basis of \(\xi\)-averages supported by \(L_0\)
(and thus also by \(P\)) is a \(c_0^{\gamma + \delta}\)-spreading model.
Since \(\beta = \gamma + \delta\), we contradict the choice of \(P\).
Therefore, this claim holds as well.

Summarizing, we showed that for every \(\gamma < \xi + \beta\) and all \(M \in [P]\)
there exists \(L \in [M]\) such that no block basis of \(\gamma\)-averages
supported by \(L\) is a \(c_0^\delta\)-spreading model, where
\(\gamma + \delta = \xi + \beta\). But this means \(P \in [N]\) is
\((\xi + \beta)\)-large, contradicting the definition of \(\xi\).
The proof of the proposition is now complete.
\end{proof}
In the next part of this section we give the proof of Theorem
\ref{mT2}. We shall need a few technical lemmas.
\begin{Lem} \label{mL1}
Suppose that \(N \in [\mathbb{N}]\) is \(\alpha\)-nice (see
Definition \ref{mD2}). Then for every \(P \in [N]\), every \(\beta
< \alpha\), every \(p \in \mathbb{N}\) and all \(\epsilon > 0\),
there exists \(M \in [P]\) such that every \(\alpha\)-average
supported by \(M\) is \((\beta, p, \epsilon)\)-large.
\end{Lem}
\begin{proof}
Define
\(\mathcal{D}= \{ L \in [P] : \, \alpha_1^L \text{ is }
(\beta, p, \epsilon)-\text{large} \}\).
Lemma \ref{L12} yields \(\mathcal{D}\) is closed in the
topology of pointwise convergence. Because \(N\)
is \(\alpha\)-nice,
we deduce that \([L] \cap \mathcal{D} \ne \emptyset\),
for all \(L \in [P]\). We infer now, from Theorem \ref{ram},
that \([M] \subset \mathcal{D}\), for some \(M \in [P]\).
Clearly, \(M\) is as desired.
\end{proof}
\begin{Lem} \label{mL2}
Suppose that \(N_1 \supset N_2 \supset \dots\) are infinite
subsets of \(\mathbb{N}\) and \(\alpha_1 < \alpha_2 < \dots\)
are countable ordinals such that
\(N_i\) is \(\alpha_i\)-nice
for all \(i \in \mathbb{N}\). Let \(N \in [\mathbb{N}]\) be
such that \(N \setminus N_i\) is finite, for all \(i \in \mathbb{N}\).
Then, \(N\) is \(\alpha\)-nice,
where \(\alpha = \lim_i \alpha_i\).
\end{Lem}
\begin{proof}
Let \(M \in [N]\), \(\beta < \alpha\), \(p \in \mathbb{N}\)
and \(\epsilon > 0\). It suffices to find an \(\alpha\)-average
\(u\) supported by \(M\) which is \((\beta, p , \epsilon)\)-large.
Choose a sequence of positive scalars \((\delta_i)\) with
\(\sum_i \delta_i < \epsilon / 6\).

Let \(k \in \mathbb{N}\) be such that \(\beta < \alpha_k\).
Since \(N_k\) is \(\alpha_k\)-nice,
we may apply Lemma \ref{mL1}, successively, to obtain infinite subsets
\(P_1 \supset P_2 \supset \dots\) of \(M \cap N_k\) such that,
for all \(i \in \mathbb{N}\), every \(\alpha_k\)-average supported by
\(P_i\) is \((\beta, p , \delta_i)\)-large.
Next choose integers \(p_1 < p_2 < \dots\) such that
\(p_i \in P_i\), for all \(i \in \mathbb{N}\), and set
\(P = (p_i)\).

We now employ Proposition \ref{Pe1} to find \(Q \in [P]\) with the property
that every \(\alpha\)-average supported by \(Q\) admits an
\((\epsilon/2, \alpha_k, \beta_k)\)-decomposition (see Definition \ref{mD3}),
where \(\alpha_k + \beta_k = \alpha\). Let \(u\) be an \(\alpha\)-average
supported by \(Q\). Write \(u = \sum_{i=1}^n \lambda_i u_i\),
where \(u_1 < \dots < u_n\) are normalized blocks, \((\lambda_i)_{i=1}^n\)
are positive scalars for which there exists \(I \subset \{1, \dots, n\}\)
satisfying
\[u_i \text{ is an } \alpha_k-\text{average for all } i \in I, \text{ while }
\bigl \| \sum_{i \in \{1, \dots, n\} \setminus I } \lambda_i u_i \bigr \|_{\ell_1} < \epsilon / 2.\]
If \(u_i = \sum_s a_s^i f_s\), for \(i \leq n\),
then, clearly,
\(\sum_{i \in \{1, \dots, n\} \setminus I } \lambda_i
\sum_s a_s^i < \epsilon / 2\).

We are going to show that \(u\) is \((\beta, p, \epsilon)\)-large.
To this end, let \(J\) be the union of less than, or equal to,
\(p\) consecutive members of \(S_\beta\) and let \(t \in K\).
Write \(I = \{i_1 < \dots, < i_m\}\). Observe that
\(u_{i_j}\) is an \(\alpha_k\)-average supported by \(P_j\) and thus by
the choice of \(P_j\),
\[ \bigl | \sum_{s \in J}
a_s^{i_j} f_s(t) \bigr | \leq \delta_j +
\sum_{s \notin J} a_s^{i_j} |f_s(t)|
, \text{ for all } j \leq m.\]
Therefore, letting \(I^c = \{1, \dots, n\} \setminus I\),
{\allowdisplaybreaks \begin{align}
\bigl | \sum_{i=1}^n \lambda_i \sum_{s \in J}
a_s^i f_s(t) \bigr | &\leq
\bigl | \sum_{i \in I^c } \lambda_i
\sum_{s \in J} a_s^i f_s(t) \bigr | +
\bigl | \sum_{i \in I} \lambda_i
\sum_{s \in J} a_s^i f_s(t) \bigr | \notag \\
&\leq \sum_{i \in I^c } \lambda_i \sum_s a_s^i +
\sum_{i \in I} \lambda_i \bigl | \sum_{s \in J} a_s^i f_s(t) \bigr | \notag \\
&\leq \epsilon /2 + \sum_{j=1}^m \lambda_{i_j} \bigl
( \delta_j + \sum_{s \notin J} a_s^{i_j} |f_s(t)| \bigr ) \notag \\
&\leq \epsilon /2 + 3 \sum_{j=1}^{|I|} \delta_j
+ \sum_{i=1}^n \lambda_i \sum_{s \notin J} a_s^i |f_s(t)| \notag \\
&\leq \epsilon + \sum_{i=1}^n \lambda_i \sum_{s \notin J} a_s^i |f_s(t)|. \notag
\end{align}}
The proof of the lemma is now complete.
\end{proof}
\begin{Lem} \label{L21}
Let \(u_1 < \dots < u_n\) be a normalized finite block basis
of \((f_i)\). Write
\(u_i = \sum_s a_s^i f_s\),
and set \(k_i = \max \mathrm{ supp } \, u_i\)
for all \(i \leq n\).
Let \(\alpha < \omega_1\) and
denote by \((\alpha_j + 1)_{j=1}^\infty\) the sequence of ordinals
associated to \(\alpha\). Let
\(\mathcal{G}\)
be a hereditary and spreading family,
and
\((\delta_i)_{i=1}^n \) be a sequence of non-negative scalars.
Suppose that \(J \in \mathcal{G}[S_\alpha]\) satisfies the following
property:
If \(2 \leq i \leq n\) is so that
\(J \cap \mathrm{ supp } \, u_i\)
is contained in the union of less than, or equal to,
\(k_{i-1}\) consecutive members of \(S_{\alpha_j}\), for
some \(j \leq k_{i-1}\) then,
\[\bigl |\sum_{s \in J} a_s^i f_s(t) \bigr | \leq \delta_i +
\sum_{s \notin J} |a_s^i| |f_s(t)|, \text{ for all }
t \in K.\]
Then for every scalar sequence \((b_i)_{i=1}^n\) and all \(t \in K\),
we have the estimate
\begin{align} \label{e1}
\bigl | \sum_{i=1}^n b_i \sum_{s \in J} a_s^i f_s(t) \bigr |
&\leq \max \biggl \{ \bigl | \sum_{i \in I} b_i u_i (t) \bigr |: \,
(u_i)_{i \in I} \text{ is } \,
\mathcal{G}^{+}-\text{admissible} \biggr \} \\
&+ \bigl ( \sum_{i=1}^n \delta_i \bigr ) \max_{i \leq n} |b_i|
+ \sum_{i=1}^n |b_i| \sum_{s \notin J} |a_s^i| |f_s(t)|. \notag
\end{align}
\end{Lem}
\begin{proof}
We may assume that
\(J \cap \cup_{i=1}^n \mathrm{ supp } \, u_i \ne \emptyset\),
or else the assertion of the lemma is trivial.
We may thus write \(J \cap
\cup_{i=1}^n \mathrm{ supp } \, u_i =
\cup_{l=1}^p J_l\), where
\(J_1 < \dots < J_p\) are non-empty members of
\(S_\alpha\) with
\(\{\min J_l : \, l \leq p \} \in \mathcal{G}\).

Define \(I_l = \{ i \leq n : \, r(u_i) \cap J_l \ne \emptyset\}\)
(where \(r(u_i)\) denotes the range of \(u_i\))
and \(i_l = \min I_l\), for all \(l \leq p\).
Put \(I= \{ i_l : \, l \leq p\}\) and let
\(I^c\) be the complement of \(I\) in \(\{1, \dots, n\}\).
Then \((u_i)_{i \in I}\) is \(\mathcal{G}^{+}\)-admissible.

Indeed, set \(L_i = \{ l \leq p : \, i_l = i \} \), for
all \(i \in I\). Observe that \(L_i\) is an interval and
that \(L_i < L_{i'}\) for all \(i < i'\) in \(I\).
Hence,
\(\min J_{\min L_i} \leq \max \mathrm{ supp } \, u_i \),
for all \(i \in I\).
Since \(\mathcal{G}\) is hereditary and spreading, we infer that
\((k_i)_{i \in I} \in \mathcal{G}\).
It follows now, by the spreading property of \(\mathcal{G}\),
that \((u_i)_{i \in I \setminus \{\min I \}}\)
is \(\mathcal{G}\)-admissible.

Next assume that \(i \in I^c \cap \cup_{l \leq p} I_l\).
Then there is a unique
\(l \leq p\) with \(i \in I_l\). Otherwise,
\(r(u_i) \cap J_l \ne \emptyset\) for at least two
distinct \(l\)'s, and so \(i \in I\).

It follows now that \(J \cap \mathrm{ supp } \, u_i
= J_l \cap \mathrm{ supp } \, u_i\),
for some \(l \leq p\). Note that \(i_l < i\) and
that \(J_l \cap r(u_{i_l}) \ne \emptyset\).
Therefore \(\min J_l \leq k_{i_l}\).
We deduce from this that \(J_l \in S_{\alpha_j + 1}\)
for some \(j \leq k_{i_l}\) and, subsequently, that
\(J_l\) is contained in the union of less than or equal to
\(k_{i_l}\) consecutive members of \(S_{\alpha_j}\),
for some \(j \leq k_{i_l}\).
The same holds for \(J \cap \mathrm{ supp } \, u_i\)
and as \(i_l < i\), we infer from
our hypothesis, that
\[\bigl | \sum_{s \in J} a_s^i f_s(t) \bigr |
 \leq \delta_i + \sum_{ s \notin J} |a_s^i| |f_s(t)|,
 \text{ for all } i \in I^c \text{ and } t \in K.\]
Now let \((b_i)_{i=1}^n\) be any scalar sequence
and let \(t \in K\). Then
\[ \sum_{i=1}^n b_i \sum_{s \in J}
a_s^i f_s(t) =
\sum_{i \in I} b_i \sum_{s \in J} a_s^i f_s(t) +
\sum_{i \in I^c} b_i \sum_{s \in J} a_s^i f_s(t).\]
Our preceding discussions yield
\begin{align} \label{e3}
\bigl | \sum_{i \in I^c} b_i \sum_{s \in J} a_s^i f_s(t) \bigr | &\leq
\sum_{i \in I^c}| b_i| \bigl | \sum_{s \in J} a_s^i f_s(t) \bigr | \\
&\leq \sum_{i \in I^c} |b_i| \bigl ( \delta_i + \sum_{ s \notin J} |a_s^i| |f_s(t)| \bigr ) \notag \\
&\leq (\max_{i \leq n} |b_i|) \sum_{i \in I^c} \delta_i +
\sum_{i \in I^c} |b_i| \sum_{ s \notin J} |a_s^i| |f_s(t)| \notag
\end{align}
and
\begin{align} \label{e4}
\bigl | \sum_{i \in I} b_i \sum_{s \in J} a_s^i f_s(t) \bigr | &=
\bigl | \sum_{i \in I} b_i \bigl ( u_i(t) -
\sum_{s \notin J} a_s^i f_s(t) \bigr ) \bigr | \\
&\leq \bigl | \sum_{i \in I} b_i u_i(t) \bigr | +
\sum_{i \in I}| b_i|
\sum_{s \notin J} |a_s^i| |f_s(t)|. \notag
\end{align}
Combining \eqref{e3} with \eqref{e4} we obtain \eqref{e1},
since \((u_i)_{i \in I}\) is \(\mathcal{G}^+\)-admissible.
\end{proof}
\begin{Lem} \label{mL3}
Suppose that \(N \in [\mathbb{N}]\) is \(\alpha\)-nice
and that there exist
\(\Gamma \in [N]\) and \(\gamma < \omega_1\) such that
no block basis of \(\alpha\)-averages supported by
\(\Gamma\) is a \(c_0^\gamma\)-spreading model.
Then there exist \(M \in [N]\) and
\(1 \leq \beta \leq \gamma\) such that
\(M\) is \((\alpha + \beta)\)-nice.
\end{Lem}
\begin{proof}
Define
\begin{align}
\beta = \min &\{ \psi < \omega_1 : \, \exists \, \Psi \in [N] \text{ such that no block basis of }
\notag \\
&\alpha-\text{averages supported by } \Psi \text{ is a } c_0^\psi-\text{spreading model} \}.
\notag
\end{align}
Our assumptions yield \(1 \leq \beta \leq \gamma\).
Choose \(M \in [N]\) such that no block basis of \(\alpha\)-averages
supported by \(M\) is a \(c_0^\beta\)-spreading model.
We are going to show that \(M\) is \((\alpha + \beta)\)-nice.
Let \(M_0 \in [M]\) and \(\tau < \alpha + \beta\).
Let \(p \in \mathbb{N}\) and \(\epsilon > 0\).
We shall exhibit an \((\alpha + \beta)\)-average supported by
\(M_0\) which is \((\tau, p , \epsilon)\)-large.
Choose a decreasing sequence of positive scalars
\((\delta_i)\) such that \(\sum_i \delta_i < \epsilon /6\).

We first consider the case \(\tau < \alpha\).
Because \(N\) is \(\alpha\)-nice,
we may apply Lemma \ref{mL1}, successively, to obtain
infinite subsets \(P_1 \supset P_2 \supset \dots\)
of \(M_0\) such that, for all \(i \in \mathbb{N}\),
every \(\alpha\)-average supported by \(P_i\) is
\((\tau, p, \delta_i)\)-large.
Choose integers \(p_1 < p_2 < \dots\)
such that \(p_i \in P_i\), for all \(i \in \mathbb{N}\),
and set \(P_0 = (p_i)\). Proposition \ref{Pe1}
now yields an \((\alpha + \beta)\)-average \(u\) supported
by \(P_0\) and admitting an \((\epsilon /2, \alpha, \beta)\)-decomposition
(see Definition \ref{mD3}).
In particular, there exist normalized blocks
\(u_1 < \dots < u_n\), positive scalars \((\lambda_i)_{i=1}^n\)
and \(I \subset \{1, \dots, n\}\) such that \(u= \sum_{i=1}^n \lambda_i u_i\),
\(u_i\) is an \(\alpha\)-average for all
\(i \in I\) and
\(\| \sum_{i \in \{1, \dots, n\} \setminus I} \lambda_i u_i \|_{\ell_1}
< \epsilon /2 \).
Let \(J\) be the union of less than, or equal to,
\(p\) consecutive members of \(S_\tau\), and let \(t \in K\).
By repeating the argument in the
last part of the proof of Lemma \ref{mL2} we conclude that
\(u\) is \((\tau, p, \epsilon)\)-large. This proves the assertion when
\(\tau < \alpha\).

Next suppose \(\alpha \leq \tau < \alpha + \beta\) and
choose \(\zeta < \beta\) with \(\tau = \alpha + \zeta\).
Recall that the definition of \(\beta\) implies that every
infinite subset of \(M_0\) supports a block basis of \(\alpha\)-averages
which is a \(c_0^\zeta\)-spreading model. Hence, thanks to Lemma \ref{L13},
there will be no loss of generality in assuming that for some positive
constant \(C\), every block basis of \(\alpha\)-averages supported by \(M_0\)
is a \(c_0^\zeta\)-spreading model with constant \(C\).
We shall further assume, because of Lemma \ref{conv}, that for every \(F \in S_{\tau}[M_0]\)
we have \(F \setminus \{\min F\} \in S_\zeta[S_\alpha]\).

Let \((\alpha_j + 1\)) be the sequence of ordinals associated to \(\alpha\).
We shall construct \(m_1 < m_2 < \dots\) in \(M_0\)
with the following property: If \(n \in \mathbb{N}\) and
\(j \leq m_n\), then every \(\alpha\)-average supported by
\(\{m_i : \, i > n\}\) is \((\alpha_j, m_n, \delta_n)\)-large.
This construction is done inductively as follows:
Choose \(m_1 \in M_0\). Apply Lemma \ref{l1} to find
\(L_1 \in [M_1]\) with \(m_1 < \min L_1\) and such that
\(S_{\alpha_j} [L_1] \subset S_{\alpha_{m_1}}\) for all
\(j \leq m_1\). We then employ Lemma \ref{mL1}, as \(N\)
is \(\alpha\)-nice, to obtain
\(M_1 \in [L_1]\) such that every \(\alpha\)-average supported by \(M_1\)
is \((\alpha_{m_1}, m_1, \delta_1)\)-large.
It follows that every \(\alpha\)-average supported by \(M_1\) is
\((\alpha_j, m_1, \delta_1)\)-large,
for all \(j \leq m_1\). Set \(m_2 = \min M_1\).

Suppose \(n \geq 2\) and that we have selected integers
\(m_1 < \dots < m_n\) in \(M_0\), and infinite subsets
\(M_1 \supset \dots \supset M_{n-1}\) of \(M_0\)
with \(m_{i+1} = \min M_i\) and such that
every \(\alpha\)-average supported by \(M_i\) is
\((\alpha_j, m_i, \delta_i)\)-large for all
\(j \leq m_i\) and \(i < n\).

We next choose, by Lemma \ref{l1}, \(L_n \in [M_{n-1}]\) with \(m_n < \min L_n\)
and such that \(S_{\alpha_j}[L_n] \subset S_{\alpha_{m_n}}\),
for all \(j \leq m_n\). Because \(N\)
is \(\alpha\)-nice,
Lemma \ref{mL1} allows us select \(M_n \in [L_n]\) such that
every \(\alpha\)-average supported by \(M_n\)
is \((\alpha_j, m_n, \delta_n)\)-large for all \(j \leq m_n\).
Set \(m_{n+1} = \min M_n\).
This completes the inductive step. Evidently, \(m_1 < m_2 < \dots\)
satisfy the required property.

We set \(P = (m_n)\). The preceding construction yields the following fact that
will be used later in the course of the proof: Suppose \(v\) is an
\(\alpha\)-average supported by \(P\) and \(\min \mathrm{ supp } \, v = m_n\),
for some \(n \geq 2\), then \(v\) is \((\alpha_j, m_{n-1}, \delta_{n-1})\)-large,
for all \(j \leq m_{n-1}\).

Recall that no block basis of \(\alpha\)-averages supported by \(P\) is
a \(c_0^\beta\)-spreading model. Let \(0 < \delta < \epsilon / (p(C+1) + 3)\)
and apply Theorem \ref{Te1} to find an \((\alpha + \beta)\)-average \(u\)
supported by \(P\), normalized blocks \(u_1 < \dots < u_n\), positive scalars
\((\lambda_i)_{i=1}^n\) and \(I \subset \{1, \dots, n\}\) such that
\(u= \sum_{i=1}^n \lambda_i u_i\),
\(u_i\) is an \(\alpha\)-average for all
\(i \in I\),
\(\| \sum_{i \in \{1, \dots, n\} \setminus I} \lambda_i u_i \|_{\ell_1}
< \delta\) and \(\max_{i \in I} \lambda_i < \delta\).
We show \(u\) is \((\tau, p, \epsilon)\)-large which will finish
the proof of the lemma.
Set
\[\mathcal{G} = \{ F \in [\mathbb{N}]^{< \infty} : \, \exists \,
F_1 < \dots < F_p \text{ in } S_\zeta^{+}, \,
F \subset \cup_{i=1}^p F_i \}.\]
\(\mathcal{G}\) is a hereditary and spreading family.

Let \(J \subset M_0\) be the union of less than, or equal to, \(p\)
consecutive members of \(S_\tau\), and let
\(t \in K\).
Our assumptions on \(M_0\) yield
\(J \in \mathcal{G}[S_\alpha]\).
Let \(\{i_1 < \dots, < i_m \}\) be an enumeration of \(I\)
and put \(m_{d_k} = \max \mathrm{ supp } \, u_{i_k}\), for all \(k \leq m\).
It has been already remarked that \(u_{i_k}\) is
\((\alpha_j, m_{d_{k-1}}, \delta_{d_{k-1}})\)-large, for all
\(2 \leq k \leq m\) and \(j \leq m_{d_{k-1}}\).
It follows that the hypotheses of Lemma \ref{L21} are fulfilled
for the block basis \(u_{i_1} < \dots < u_{i_m}\) and the given \(J \subset M_0\),
with ``\(\delta_1 \)''\(= 0\) and ``\(\delta_k\)''\( = \delta_{d_{k-1}}\)
for \(2 \leq k \leq m\). Writing \(u_i = \sum_s a_s^i f_s\), for all \(i \leq n\),
we infer from \eqref{e1} that
{\allowdisplaybreaks \begin{align}
\bigl | \sum_{i \in I} \lambda_i \sum_{s \in J} a_s^i f_s(t) \bigr |
\leq \max &\biggl \{ \bigl | \sum_{i \in E} \lambda_i u_i (t) \bigr |: \,
E \subset I, \, (u_i)_{i \in E}
\text{ is } \notag \\
&\mathcal{G}^{+}-\text{admissible} \biggr \} \notag \\
&+ \bigl (
\sum_{i=1}^\infty \delta_i \bigr ) \max_{i \in I} \lambda_i
+ \sum_{i \in I} \lambda_i \sum_{s \notin J} a_s^i |f_s(t)|. \notag
\end{align}}
Note that when \((u_i)_{i \in E}\) is \(\mathcal{G}^{+}\)-admissible, we have
\[\bigl \| \sum_{i \in E} \lambda_i u_i \bigr \| \leq \bigl
(p(C + 1) + 1 \bigr ) \max_{i \in E} \lambda_i < \bigl (p(C + 1) + 1 \bigr ) \delta.\]
Hence,
\[\bigl | \sum_{i \in I} \lambda_i \sum_{s \in J} a_s^i f_s(t) \bigr |
< \bigl (p(C + 1) + 2 \bigr ) \delta
+ \sum_{i \in I} \lambda_i \sum_{s \notin J} a_s^i |f_s(t)|.\]
Next, put \(I^c = \{1, \dots, n\} \setminus I\).
Then,
\[\sum_{i \in I^c} \lambda_i \sum_s  a_s^i
< \delta, \text{ as }
\bigl \| \sum_{i \in I^c} \lambda_i u_i \bigr \|_{\ell_1}
< \delta.\]
Combining the preceding estimates we conclude
\begin{align}
\bigl | \sum_{i=1}^n \lambda_i \sum_{s \in J}
a_s^i f_s(t) \bigr | &\leq
\sum_{i \in I^c } \lambda_i
\sum_s a_s^i +
\bigl | \sum_{i \in I} \lambda_i
\sum_{s \in J}
a_s^i f_s(t) \bigr | \notag \\
&< \delta  + \bigl (p(C + 1) + 2 \bigr ) \delta  +
\sum_{i \in I} \lambda_i \sum_{s \notin J} a_s^i |f_s(t)| \notag \\
&< \epsilon + \sum_{i=1}^n \lambda_i \sum_{s \notin J} a_s^i |f_s(t)|.
\notag
\end{align}
Therefore,
\(u\) is \((\tau, p, \epsilon)\)-large.
This completes the proof.
\end{proof}
We are now ready for the
\begin{proof}[Proof of Theorem \ref{mT2}]
We claim that every infinite subset of \(N\)
contains a further infinite subset which is \(\alpha\)-nice.
If this claim holds, then
evidently, \(N\) is itself \(\alpha\)-nice.
So suppose on the contrary, that the claim is false and choose
\(N_0 \in [N]\) having no infinite subset which is
\(\alpha\)-nice.
We now claim that there exist \(1 \leq \beta_1 < \alpha\)
and \(N_1 \in [N_0]\) which is \(\beta_1\)-nice.
Indeed, define
\begin{align}
\beta_1 = \min &\{ \zeta < \omega_1 : \, \exists \, M \in [N_0] \text{ such that no block basis of }
\notag \\
&0-\text{averages supported by } M \text{ is a } c_0^\zeta-\text{spreading model} \}.
\notag
\end{align}
Since \(N\) is \(\alpha\)-large, \(\alpha\) belongs to the set and so
\(1 \leq \beta_1 \leq \alpha\).
Choose \(N_1 \in [N_0]\) such that no block basis of \(0\)-averages
supported by \(N_1\) is a \(c_0^{\beta_1}\)-spreading model.
We show \(N_1\) is \(\beta_1\)-nice.
Because \(N_0\) is assumed to contain no infinite subset which is
\(\alpha\)-nice, we shall also obtain \(\beta_1 < \alpha\).

Let \(M \in [N_1]\), \(\beta < \beta_1\), \(p \in \mathbb{N}\) and \(\epsilon > 0\).
We shall find a \(\beta_1\)-average supported by \(M\) which is \((\beta, p, \epsilon)\)-large.
Since \(\beta < \beta_1\), there exist \(M_1 \in [M]\) and a constant \(C > 0\) such that
the block basis \((f_m)_{m \in M_1}\) is a \(c_0^\beta\)-spreading model with constant \(C > 0\).
Let \(0 < \delta < \epsilon / (pC)\). Since no block basis of \(0\)-averages supported by \(M_1\)
is a \(c_0^{\beta_1}\)-spreading model, Theorem \ref{Te1} yields a
\(\beta_1\)-average \(u\), supported by \(M_1\), positive scalars \((\lambda_i)_{i \in F}\)
(where \(F = \mathrm{ supp } \, u\)) and \(I \subset F\) with \(I \in S_{\beta_1}\),
such that
\[u = \sum_{i \in F} \lambda_i f_i, \, \max_{i \in I} \lambda_i < \delta, \, \text{ and }
\sum_{i \in F \setminus I} \lambda_i < \delta.\]
Let \(t \in K\) and let \(J\) be the union of less than,
or equal to, \(p\) consecutive members of \(S_\beta\).
It follows that
\begin{align}
\bigl | \sum_{ i \in J \cap F} \lambda_i f_i (t) \bigr | &\leq
\bigl \|\sum_{ i \in J \cap F} \lambda_i f_i \bigr \| \notag \\
&\leq pC \max_{i \in F} \lambda_i < pC \delta < \epsilon. \notag
\end{align}
Thus, \(u\) is a \(\beta_1\)-average, \((\beta, p, \epsilon)\)-large, and so
\(N_1\) is \(\beta_1\)-nice, as claimed.

We shall now construct, by transfinite induction on \(1 \leq \tau < \omega_1\),
families \(\{N_\tau\}_{1 \leq \tau < \omega_1} \subset [N_0] \)
and \(\{\beta_\tau\}_{1 \leq \tau < \omega_1} \subset [1, \alpha)\)
with the following properties:
\begin{enumerate}
\item \(N_{\tau_2} \setminus N_{\tau_1} \) is finite, for all
\(1 \leq \tau_1 < \tau_2 < \omega_1\).
\item \(N_\tau\) is \(\beta_\tau\)-nice,
for all \(1 \leq \tau < \omega_1\).
\item \(\beta_{\tau_1} < \beta_{\tau_2}\), for all
\(1 \leq \tau_1 < \tau_2 < \omega_1\).
\end{enumerate}
Of course, \((3)\) is absurd since \(\alpha < \omega_1\). Hence,
our assumption that \(N_0\) contained no infinite subset which is
\(\alpha\)-nice, was false. The proof of the theorem will be
completed, once we give the construction of the above described
families, satisfying conditions \((1)\)-\((3)\). \(N_1\) and
\(\beta_1\) have been already constructed. Suppose that \(1 <
\tau_0 < \omega_1\) and that \(\{N_\tau\}_{1 \leq \tau < \tau_0}
\subset [N_0] \), \(\{\beta_\tau\}_{1 \leq \tau < \tau_0} \subset
[1, \alpha)\) have been constructed fulfilling properties
\((1)\)-\((3)\), above, with \(\omega_1\) being replaced by
\(\tau_0\).

Assume first that \(\tau_0\) is a successor ordinal, say
\(\tau_0 = \tau_1 + 1\). We know by the inductive construction, that
\(N_{\tau_1}\) is \(\beta_{\tau_1}\)-nice.
By assumption, \(N\) is \(\alpha\)-large. Since
\(\beta_{\tau_1} < \alpha\), there exists \(\Gamma \in [N_{\tau_1}]\)
such that no block basis of \(\beta_{\tau_1}\)-averages supported by \(\Gamma\)
is a \(c_0^{\eta_{\tau_1}}\)-spreading model, where
\(\beta_{\tau_1} + \eta_{\tau_1} = \alpha\).
Lemma \ref{mL3} now implies the existence of
\(N_{\tau_0} \in [N_{\tau_1}]\) and
\(1 \leq \zeta_{\tau_1} \leq \eta_{\tau_1}\)
such that \(N_{\tau_0}\) is \((\beta_{\tau_1} + \zeta_{\tau_1})\)-nice.
Set \(\beta_{\tau_0} = \beta_{\tau_1} + \zeta_{\tau_1}\).
Necessarily, \( \beta_{\tau_0} < \alpha\), by the choice of \(N_0\).
It is easy to see that the families
\(\{N_\tau\}_{1 \leq \tau < \tau_0 + 1}\) and
\(\{\beta_\tau\}_{1 \leq \tau < \tau_0 + 1}\)
satisfy conditions \((1)\)-\((3)\), above, with
\(\omega_1\) being replaced by \(\tau_0 +1\).

Next assume that \(\tau_0\) is a limit ordinal and
choose a strictly increasing sequence of ordinals
\(\tau_1 < \tau_2 < \dots\)
such that \(\tau_0 = \lim_n \tau_n\).
By the inductive construction we have that
\(\beta_{\tau_1} < \beta_{\tau_2} < \dots\) and
thus we may define the limit ordinal
\(\beta_{\tau_0} = \lim_n \beta_{\tau_n}\).
In addition to this,
\(N_{\tau_n} \setminus N_{\tau_m}\) is finite
for all integers \(m < n\). We deduce from the above, that
\(\cap_{i=1}^k N_{\tau_i}\) is
\(\beta_{\tau_k}\)-nice,
for all \(k \in \mathbb{N}\). Finally, choose
\(N_{\tau_0} \in [N_0]\) such that
\(N_{\tau_0} \setminus \cap_{i=1}^k N_{\tau_i}\)
is finite, for all \(k \in \mathbb{N}\).
We infer from Lemma \ref{mL2},
that \(N_{\tau_0}\) is \(\beta_{\tau_0}\)-nice.
It is easily verified now, that
the families
\(\{N_\tau\}_{1 \leq \tau < \tau_0 + 1}\) and
\(\{\beta_\tau\}_{1 \leq \tau < \tau_0 + 1}\)
satisfy conditions \((1)\)-\((3)\), above, with
\(\omega_1\) being replaced by \(\tau_0 +1\).
This completes the inductive step and the proof
of the theorem.
\end{proof}
\begin{proof}[Proof of Corollary \ref{elton}]
Assume without loss of generality, that
\((f_n)\) has no subsequence equivalent to the unit vector basis
of \(c_0\). By the Kunen-Martin boundedness principle
(see \cite{d}, \cite{ke}), we may
choose an ordinal \(1 \leq \gamma < \omega_1\) such that no subsequence
of \((f_n)\) is a \(c_0^\gamma\)-spreading model.
Set \(K_m = \{ t \in K : \, \sum_n |f_n(t)| \leq m \}\),
for all \(m \in \mathbb{N}\). Clearly, \((K_m)\) is
an increasing sequence of closed subsets of \(K\)
and \(K = \cup_m K_m\). We claim that
for every \(m \in \mathbb{N}\), every \(N \in [\mathbb{N}]\),
and all \(\epsilon > 0\),
there exists a \(\gamma\)-average \(u\) of \((f_n)\)
supported by \(N\) and such that \(|u|(t) < \epsilon\),
for all \(t \in K_m\)
(if \(u = \sum_i a_i f_i\), we define
\(|u|(x) = \sum_i |a_i| |f_i(x)|\), for all \(x \in K\)).

To see this, let \(0 < \delta < \epsilon /m\).
Since no subsequence of \((f_n)\) is a \(c_0^\gamma\)-spreading model,
Theorem \ref{Te1} allows us choose a \(\gamma\)-average \(u\)
of \((f_n)\), supported by \(N\)
and such that there exist non-negative scalars
\((\lambda_i)_{i=1}^p\) and
\(I \subset \{1, \dots, p\} \)
satisfying the following: \((1)\)
\(u = \sum_{i=1}^p \lambda_i f_i\) and \(\max_{i \in I} \lambda_i < \delta\).
\((2)\) \((f_i)_{i \in I}\) is \(S_\gamma\)-admissible (i.e. \(I \in S_\gamma\)) and
\(\sum_{i \in \{1, \dots , p\} \setminus I} \lambda_i < \delta\).
It is easy to check now that for every  \(t \in K_m\) we have
\(|u|(t) < \epsilon\) and thus our claim holds.

Now let \((\epsilon_n)\) be a summable sequence of positive scalars
and \(N \in [\mathbb{N}]\).
Successive applications of the previous claim yield a block basis
\(v_1 < v_2 < \dots\) of \(\gamma\)-averages of \((f_n)\), supported by \(N\)
and satisfying \(|v_n| (t) < \epsilon_n\) for every \(t \in K_n\) and
all \(n \in \mathbb{N}\). It follows that for all \(t \in K\)
the set \(\{n \in \mathbb{N} : \, |v_n(t)| \geq \epsilon_n \}\)
is a subset of \(\{1, \dots, q_t\}\), where \(q_t\) is the least
\(m \in \mathbb{N}\) such that \(t \in K_m\).
We deduce from Theorem \ref{mT}, that there exist
\(\beta < \omega_1\) and a block basis of \(\beta\)-averages
of \((v_n)\), equivalent to the unit vector basis of \(c_0\).

In order to get a block basis of averages of \((f_i)\) equivalent
to the unit vector basis of \(c_0\), one needs a somewhat more
demanding argument which goes as follows.
Choose a countable limit ordinal \(\alpha\) with \(\gamma < \alpha \)
and let \((\alpha_j + 1)_{j=1}^\infty\) be the sequence of ordinals
associated to \(\alpha\). Let \(N \in [\mathbb{N}]\) and choose
\(n \in N\) with \(n \geq 2\), such that
\(\gamma < \alpha_n\). Let \(m \in \mathbb{N}\).
Since no subsequence of \((f_i)\) is a \(c_0^{\alpha_n}\)-spreading model,
our preceding argument allows us choose an \(\alpha_n\)-average \(v\)
of \((f_i)\), supported by \(\{i \in N : \, n < i \}\),
and such that \(|v|(t) < 1/(2n)\), for all
\(t \in K_m\). Set
\(u = \bigl ( (1/n) f_n + v \bigr ) /
\| (1/n) f_n + v\|\).
Clearly, \(u\) is an \(\alpha\)-average of \((f_i)\)
supported by \(N\) and satisfying
\(|u|(t) < 3/n\), for all
\(t \in K_m\). Note that \(n = \min \mathrm{ supp } \, u\).

Summarizing, given \(N \in [\mathbb{N}]\) we can select a block basis
\(u_1 < u_2 < \dots\) of \(\alpha\)-averages of
\((f_i)\) supported by \(N\) and satisfying
\(|u_n|(t) < 3/m_n\), for all \(t \in K_n\)
and \(n \in \mathbb{N}\). In the above, we have
let \(m_n = \min \mathrm{ supp } \, u_n\),
for all \(n \in \mathbb{N}\). It follows that for all \(n \in \mathbb{N}\),
if \(t \in K_n\) and
\(|u_i|(t) \geq 3/m_i\), then \(i < n\).
Given \(L \in [\mathbb{N}]\), set \(l_n = \min \mathrm{ supp } \, \alpha_n^L\),
for all \(n \in \mathbb{N}\).
We now define
\[\mathcal{D} = \{ L \in [\mathbb{N}] : \, \forall \, n \in \mathbb{N}, \,
\forall \, t \in K_n,
\text{ if } |\alpha_i^L|(t) \geq 3/l_i, \text{ then }
i < n \}. \]
\(\mathcal{D}\) is closed in the topology of pointwise convergence,
thanks to Lemma \ref{L12}. Our preceding discussion and Lemma \ref{L12}, show that
every \(N \in [\mathbb{N}]\) contains some \(L \in \mathcal{D}\) as a
subset. We infer from Theorem \ref{ram}, that
\([N] \subset \mathcal{D}\) for some \(N \in [\mathbb{N}]\).

Next, let \(\mathcal{T}_0\) be the collection of those
finite subsets \(E\) of \(N\) that can be written in the form
\(E = \cup_{i=1}^m \mathrm{ supp } \, \alpha_i^L\),
for some \(L \in [N]\) (depending on \(E\))
for which there exists some \(t \in K\) (depending
on \(E\) and \(L\)) such that
\(|\alpha_i^L|(t) \geq 3/l_i\), for all \(i \leq m\).

Let \(\mathcal{T}\) be the collection of all initial segments of
elements of \(\mathcal{T}_0\). We claim that
\(\mathcal{T}\) is compact in the topology of pointwise convergence.
Indeed, were this false, there would exist \(M \in [N]\),
\(M=(m_i)\), such that \(\{m_1, \dots, m_n\} \in \mathcal{T}\),
for all \(n \in \mathbb{N}\).
Let \(n \in \mathbb{N}\). It follows that
\(\cup_{i=1}^n \mathrm{ supp } \, \alpha_i^M
\in \mathcal{T}\). Hence, there exist
\(L_n \in [N]\), \(k_n \in \mathbb{N}\) and
\(t_n \in K\) such that
\(\cup_{i=1}^n \mathrm{ supp } \, \alpha_i^M \)
is an initial segment of
\(\cup_{i=1}^{k_n} \mathrm{ supp } \, \alpha_i^{L_n}\)
and \(|\alpha_i^{L_n}|(t_n) \geq 3/d_i\), for
all \(i \leq k_n\), where
\(d_i = \min \mathrm{ supp } \, \alpha_i^{L_n}\),
for all \(i \in \mathbb{N}\).
We now deduce from Lemma \ref{L12}, that
\(n \leq k_n\) and that \(\alpha_i^M = \alpha_i^{L_n}\),
for all \(i \leq n\). Therefore,
\(|\alpha_i^M|(t_n) \geq 3/m_i\), for all \(i \leq n\),
where
\(m_i = \min \mathrm{ supp } \, \alpha_i^M\),
for all \(i \in \mathbb{N}\).
The compactness of \(K\) now implies that
there is some \(t \in K\) satisfying
\(|\alpha_i^M|(t) \geq 3/m_i\), for all \(i \in \mathbb{N}\).
This is a contradiction, as \(M \in \mathcal{D}\).
Thus, our claim holds and so \(\mathcal{T}\) is indeed compact.

We next apply a result from \cite{o1} to obtain
\(P \in [N]\) such that \(\mathcal{T}[P]\)
is a hereditary and compact family.
The result in \cite{g} now yields
\(Q \in [P]\) and a countable ordinal
\(\eta > \alpha\), such that
\(\mathcal{T}[Q] \subset S_\eta\).
It follows that for every \(L \in [Q]\)
and all \(n \in \mathbb{N}\) such that
there exists some \(t \in K\) satisfying
\(|\alpha_i^L|(t) \geq 3/l_i\), for all
\(i \leq n\), we have
\(\cup_{i=1}^n \mathrm{ supp } \, \alpha_i^L \in S_\eta\).

We now claim that \(\xi^Q \leq \eta\) (see Definition \ref{mD1}).
If this is not the case, we may choose
\(R \in [Q]\), \(R = (r_i)\), which is \(\zeta\)-large, for some countable
ordinal \(\zeta\) with
\(\eta < \zeta\). Let \(\epsilon > 0\).
We shall assume, as we clearly may, that
\(\sum_i (1/r_i) < \epsilon\).
Since \(\alpha < \eta\), we
may choose an ordinal \(\beta\) with
\(\alpha + \beta = \zeta\).
By passing to an infinite subset of
\(R\), if necessary, we may assume without loss
of generality, thanks to Proposition \ref{Pe1},
that every \(\zeta\)-average of \((f_i)\)
supported by \(R\) admits an
\((\epsilon, \alpha, \beta)\)-decomposition.

Because \(R\) is \(\zeta\)-large, it is also
\(\zeta\)-nice, by Theorem \ref{mT2}.
We may thus select a \(\zeta\)-average
\(u\) of \((f_i)\),
supported by \(R\), which is \((\eta, 1, \epsilon)\)-large.
We infer from Proposition \ref{Pe1} that
there exist normalized blocks
\(u_1 < \dots < u_n\), positive scalars \((\lambda_i)_{i=1}^n\)
and \(I \subset \{1, \dots, n\}\) such that
\(u = \sum_{i=1}^n \lambda_i u_i\) and
\(u_i\) is an \( \alpha\)-average for all
\(i \in I\), while
\(\| \sum_{i \notin I} \lambda_i u_i \|_{\ell_1}
< \epsilon\).

Now let \(t \in K\) and define
\(H = \{i \in I : \, |u_i|(t) \geq 3 / q_i \}\),
where \(q_i = \min \mathrm{ supp } \, u_i\),
for all \(i \in I\). Let \(\{i_1 < \dots, < i_k\}\)
be an enumeration of \(H\).
Lemma \ref{L12} yields some \(L \in [R]\) such that
\(u_{i_j} = \alpha_j^L\), for all \(j \leq k\).
Set \(J = \cup_{i \in H} \mathrm{ supp } \, u_i\).
Since \(L \in [Q]\), it follows that
\(J \in S_\eta\). Writing \(u_i = \sum_s a_s^i f_s\), for all \(i \leq n\),
we conclude, as \(u\) is \((\eta, 1, \epsilon)\)-large, that
\[ \bigl |\sum_{i \in H} \lambda_i \sum_s  a_s^i f_s(t) \bigr | \leq
\epsilon + \sum_{i \notin H} \lambda_i \sum_s a_s^i |f_s(t)|.\]
We now have the estimates
{\allowdisplaybreaks \begin{align}
|u(t)| &= \bigl | \sum_{i \in H }\lambda_i \sum_s
a_s^i f_s(t) + \sum_{i \notin H} \lambda_i \sum_s
a_s^i f_s(t) \bigr | \notag \\
&\leq \bigl |\sum_{i \in H }\lambda_i \sum_s
a_s^i f_s(t) \bigr | +
\bigl |\sum_{i \notin H }\lambda_i \sum_s
a_s^i f_s(t) \bigr | \notag \\
&\leq \epsilon +  2 \sum_{i \notin H }\lambda_i \sum_s
a_s^i |f_s(t)| \notag \\
&\leq \epsilon + 2 \sum_{i \in I \setminus H }\lambda_i \sum_s
a_s^i |f_s(t)| +
2 \sum_{i \notin I}\lambda_i \sum_s
a_s^i |f_s(t)| \notag \\
&\leq \epsilon + 6 \sum_{i \in I \setminus H} |u_i|(t)
+ 2 \bigl \| \sum_{ i \notin I} \lambda_i u_i \bigr \|_{\ell_1}
\notag \\
&< \epsilon + 18 \sum_{i \in I \setminus H} (1/q_i) + 2 \epsilon
\notag \\
&< 21 \epsilon. \notag
\end{align}}
Since \(\|u\| = 1\), we reach a contradiction for
\(\epsilon\) small enough.
Therefore, \(\xi^Q \leq \eta\). Proposition
\ref{mP} now yields a block basis of
\(\xi\)-averages of \((f_i)\), for some \(\xi \leq \eta\),
equivalent to the unit vector basis of \(c_0\).
\end{proof}

\end{document}